\documentclass[11pt]{article}

\usepackage{geometry}                
\usepackage{graphicx}
\usepackage{amsmath,amssymb,amsfonts,amsthm,bm,mathtools}
\usepackage{tikz-cd} 
\usepackage{cite}
\usepackage{authblk} 
\usepackage{hyperref} 
\usepackage{upgreek}
\usepackage{bbm}
\usepackage{pifont}
\usepackage{subfigure}
\usepackage{subcaption}
\usepackage{stmaryrd} 
\SetSymbolFont{stmry}{bold}{U}{stmry}{m}{n} 

\hypersetup{%
  colorlinks=true,
  linkcolor=blue,
  linkbordercolor=red,
  citecolor=red,
}

\makeatletter
\DeclareRobustCommand{\pdot}{\mathbin{\mathpalette\pdot@\relax}}
\newcommand{\pdot@}[2]{%
  \ooalign{%
    $\m@th#1\circ$\cr
    \hidewidth$\m@th#1\cdot$\hidewidth\cr
  }%
}
\makeatother

\newtheorem{theorem}{Theorem}
\newtheorem{prop}{Proposition}

\newcommand{\bmsf}[1]{\boldsymbol{\mathsf{#1}}}  
\newcommand{\msf}[1]{\mathsf{#1}}

\title{A mimetic discretization of Westervelt's equation}
\author[1]{William Barham}
\author[1,2]{Philip J. Morrison}
\affil[1]{Institute for Fusion Studies, The University of Texas at Austin}
\affil[2]{Department of Physics, The University of Texas at Austin}
\date{\today} 

\begin{document}

\maketitle

\section*{Abstract}

A broad class of nonlinear acoustic wave models possess a Hamiltonian structure in their dissipation-free limit and a gradient flow structure for their dissipative dynamics. This structure may be exploited to design numerical methods which preserve the Hamiltonian structure in the dissipation-free limit, and which achieve the correct dissipation rate in the spatially-discrete dissipative dynamics. Moreover, by using spatial discretizations which preserve the de Rham cohomology, the non-evolving involution constraint for the vorticity may be exactly satisfied for all of time. Numerical examples are given using a mimetic finite difference spatial discretization.

\tableofcontents

\section{Introduction}

Nonlinear acoustic wave models may be expressed as Hamiltonian field theories in their dissipation-free limit. This quality provides a powerful arsenal of theoretical tools to study these models. The present work leverages this Hamiltonian structure to design structure-preserving numerical methods which automatically preserve known properties of the continuous model. In particular, this work studies structure-preserving discretizations of a nonlinear acoustic model known as Westervelt's equation \cite{hamilton1998nonlinear, NikolićSaid-Houari2020}. 

The presentation of the discretization framework is left general enough to accommodate structure-preserving discretizations based on both Galerkin (e.g.\ finite element) and collocation (e.g. \ finite difference) approaches. The terminology and notation follows that used in \cite{bochev_and_hyman}. As such, the general approach taken herein to accomplish the spatial discretization of Westervelt's equation may be applied to derive a finite element exterior calculus (FEEC) discretization \cite{ArnoldDouglasN2010Feec}, a mimetic spectral element method \cite{kreeft2011mimetic}, a structure-preserving Fourier spectral discretization \cite{10.1007/s10915-022-01781-3}, or a mimetic finite difference discretization \cite{bochev_and_hyman}. This work utilizes a lowest-order mimetic finite difference method and studies the problem on a periodic domain. This is done not because the method is limited, but to present a minimally simple implementation in this initial work. Higher-order methods and more complex geometries are deferred to a future work. The use of a collocation method rather than a Galerkin method is convenient due to the fact that a nonlinear system must be solved at each time-step: a lowest-order finite difference approximation yields an analytically solvable nonlinear solve. The method described herein generalizes to other nonlinear acoustic models such as Kuznetsov's equation \cite{NikolićSaid-Houari2020}, and is similar to a previously derived finite element discretization for nonlinear optical models \cite{Barham_dissertation}. The numerical examples in this work emphasize the convergence of the method, and its preservation of properties of the continuous model, i.e.\ conservation of vorticity and the known energy dissipation rate. 

The work is structured as follows. Section \ref{sec:cts_ham_theory} describes the continuous Hamiltonian structure of nonlinear acoustic models, in particular, Westervelt's equation. Section \ref{sec:disc_ham_theory} describes a general approach to the structure-preserving discretization of nonlinear acoustic models. Section \ref{sec:disc_westervelt} specializes the discussion to describe a structure-preserving discretization of Westervelt's equation. Section \ref{sec:numerical_examples} provides numerical examples in one and two spatial dimensions. Appendix \ref{appendix:mimetic_fd} describes the particular mimetic finite difference method used in this work in greater detail. Appendix \ref{appendix:nonlinear_stability} provides further explanation for the stability of the time-stepping scheme for Westervelt's equation. 

\section{The continuous Hamiltonian theory of nonlinear acoustics} \label{sec:cts_ham_theory}

Many nonlinear acoustic models may be expressed as Hamiltonian systems in the dissipation-free limit with the dissipative part of the dynamics expressed as a gradient flow generated by the Hamiltonian. This section describes this general modeling framework and how Westervelt's equation in particular fits into this framework. 

\subsection{Hamiltonian field theories and dissipative dynamics}

Hamiltonian systems are those equations whose evolution is specified by a manifold $\mathcal{M}$, a Poisson bracket, $\{ \cdot, \cdot \}: C^\infty(\mathcal{M}) \times C^\infty(\mathcal{M}) \to C^\infty(\mathcal{M})$, and an energy functional, the Hamiltonian, $H: C^\infty(\mathcal{M}) \to \mathbb{R}$. The evolution for any $F \in C^\infty(\mathcal{M})$ is given by
\begin{equation}
    \dot{F} = \{ F, H \} \,.
\end{equation}
Such a construction may be used to model ordinary or partial differential equations and has been used to describe the non-dissipative limit of nearly every branch of classical physics. A general overview of Hamiltonian field theories may be found in \cite{pjm82, pjm98}. 

One calls $(\mathcal{M}, \{ \cdot, \cdot \})$ a Poisson manifold if $\mathcal{M}$ is a smooth manifold and $\{\cdot, \cdot \}: C^\infty(\mathcal{M}) \times C^\infty(\mathcal{M}) \to C^\infty(\mathcal{M})$ is a Poisson bracket. A smooth functional of the manifold $\mathcal{M}$ will often be referred to as an observable. A Poisson bracket is bilinear map on the vector space of smooth functionals of $\mathcal{M}$ such that
\begin{equation}
\begin{aligned}
    &\{F, G\} = - \{G, F\} \\
    &\{ F, \{G, H\} \} + \{G, \{H, F \} \} + \{ H, \{F, G\} \} = 0 \\
    &\{F, G H \} = G \{F, H \} + \{F, G \} H \,,
\end{aligned}
\end{equation}
$\forall F,G,H \in C^\infty(\mathcal{M})$. The first two properties make the space of smooth functions on $\mathcal{M}$ a Lie algebra with the Poisson bracket as its Lie bracket, while the final identity is the Leibniz rule which makes the Poisson bracket a derivation on $C^\infty(\mathcal{M})$. An interesting feature of noncanonical Hamiltonian mechanics is that the Poisson bracket may possess degeneracy. This generates constants of motion because
\begin{equation}
    \{F, C \} = 0 \, \quad \forall F \in C^\infty(\mathcal{M}) \implies \dot{C} = \{H, C\} = 0 \,.
\end{equation}
Such invariants are called Casimir invariants. 

One frequently couples dissipative dynamics to Hamiltonian dynamics using the metriplectic formalism \cite{MORRISON1984423, morrison2023inclusive}, however the Poisson structure of the models considered in this work is simple with no suitable entropy among their Casimir invariants. Instead, this work defines a symmetric positive definite bracket, $(\cdot, \cdot)$, such that the flow is generated by
\begin{equation}
    \dot{F} = \{F, H \} - (F, H) \,, \quad \forall F \,.
\end{equation}
A symmetric bracket of this form was previously developed to model dissipation in resistive reduced MHD \cite{morrison_hazeltine_84}. This breaks the energy conservation law yielding a model which monotonically dissipates energy:
\begin{equation}
    \dot{H} = - (H, H) \leq 0 \,.
\end{equation}
Such a system is typically interpreted as losing energy to unresolved scales. 

\subsection{General nonlinear acoustic models}

A general nonlinear acoustic models take the form
\begin{equation} \label{eq:general_nl_acoustic}
	\begin{cases}
		\rho_t + \nabla \cdot \bm{m} = 0 & \\
		\bm{v}_t + \nabla p = 0 & \,,
	\end{cases}
\end{equation}
where there exists a functional $K = K[p, \bm{v}]$ such that
\begin{equation}
	\bm{m} = \bm{v} + \frac{\delta K}{\delta \bm{v}}
	\quad \text{and} \quad
	\rho = p - \frac{\delta K}{\delta p} \,.
\end{equation}
That is, the momentum and density variables are a near identity transformation from the velocity and pressure variables. Both Westervelt's and Kuznetsov's equations, standard models in nonlinear acoustics, may be shown to fit into this general framework \cite{Barham_dissertation}. 

Define the following Hamiltonian:
\begin{equation}
	H[\rho, \bm{v}] = K - \int_\Omega \frac{\delta K}{\delta p} p \, \mathsf{d}^3 \bm{x} + \frac{1}{2} \int_\Omega \left( p^2 + | \bm{v} |^2 \right) \mathsf{d}^3 \bm{x} \,,
\end{equation}
where $\rho = \rho(p, \bm{v})$ is defined implicitly. One may show that
\begin{equation}
	\frac{\delta H}{\delta \rho} = p
	\quad \text{and} \quad
	\frac{\delta H}{\delta \bm{v}} = \bm{m} \,.
\end{equation} 
The Poisson bracket is simply the acoustic wave bracket which appears as a piece of the Poisson bracket of the full Navier-Stokes Poisson bracket \cite{PhysRevLett.45.790, pjm98}:
\begin{equation} \label{eq:poisson_bracket_general_acoustic}
	\{F, G\} = \int_\Omega \left( \frac{\delta F}{\delta \bm{v}} \cdot \nabla \frac{\delta G}{\delta \rho} - \frac{\delta G}{\delta \bm{v}} \cdot \nabla \frac{\delta F}{\delta \rho} \right) \mathsf{d}^3 \bm{x} \,.
\end{equation}
The equation \ref{eq:general_nl_acoustic} is generated by $\dot{F} = \{F, \overline{H} \}$.

Dissipation is frequently added to the model by defining a symmetric bracket
\begin{equation} \label{eq:acoustic_symmetric_bracket}
    \left( F, G \right) = - b \int_\Omega \nabla \frac{\delta F}{\delta \rho} \cdot \nabla \frac{\delta G}{\delta \rho} \mathsf{d}^3 \bm{x} \,,
\end{equation}
where $b$ is a parameter controlling the rate of dissipation. The flow is then given by $\dot{F} = \{F, H\} + (F, H)$. One has that energy is dissipated at the rate $\dot{H} = (H,H) = - b \| \nabla p \|_{L^2(\Omega)}^2$.

\subsection{An abstract formulation of the Hamiltonian structure of nonlinear acoustics}

It is helpful to describe the Hamiltonian formalism of nonlinear acoustics more abstractly at this juncture. This will be seen to facilitate the structure-preserving discretization of the model. To begin, define Hilbert spaces $V^\ell \subseteq L^2 \Lambda^\ell(\Omega)$, the Hilbert space of $L^2$-integrable differential $\ell$-forms. Moreover, define $V^*_\ell$, the dual space of functionals on the differential $\ell$-forms. Finally, let $\mathcal{H}_\ell : V^\ell \to V^*_\ell$ be the duality map (guaranteed by the Riesz representation theorem). 

Let $p \in V^0$ and $\bm{v} \in V^1$. Define $K: V^0 \times V^1 \to \mathbb{R}$, and define
\begin{equation}
    \bm{m} = \mathcal{H}_1 \bm{v} + \frac{\delta K}{\delta \bm{v}} \in V^*_1 \,,
    \quad \text{and} \quad
    \rho = \mathcal{H}_0 p - \frac{\delta K}{\delta p} \in V^*_0 \,.
\end{equation}
The Hamiltonian is defined by
\begin{equation}
    H[\rho, \bm{v}] 
    = K 
    - \left\langle \frac{\delta K}{\delta p}, p \right\rangle
    + \frac{1}{2} 
    \left( 
    \left\langle \mathcal{H}_0 p, p \right\rangle 
    +
    \left\langle \mathcal{H}_1 \bm{v}, \bm{v} \right\rangle
    \right) \,,
\end{equation}
where $p = p(\rho, \bm{v})$ is defined implicitly through a coordinate change, and the angle brackets indicate the evaluation pairing. It is possible to show that
\begin{equation}
	\frac{\delta H}{\delta \rho} = p
	\quad \text{and} \quad
	\frac{\delta H}{\delta \bm{v}} = \bm{m} \,.
\end{equation} 
The Poisson bracket is given by
\begin{equation} \label{eq:absract_pb}
    \{F, G\} 
    = 
    \left\langle \frac{\delta F}{\delta \bm{v}}, \mathsf{d}_0 \frac{\delta G}{\delta \rho} \right\rangle 
    - 
    \left\langle \frac{\delta G}{\delta \bm{v}}, \mathsf{d}_0 \frac{\delta F}{\delta \rho} \right\rangle \,,
\end{equation}
where $\mathsf{d}_\ell: V^\ell \to V^{\ell+1}$ is the exterior derivative. Finally, the dissipative bracket is given by
\begin{equation}
    (F,G)
    =
    - 
    b 
    \left\langle \mathcal{H}_1 \mathsf{d}_0 \frac{\delta G}{\delta \rho}, \mathsf{d}_0 \frac{\delta F}{\delta \rho} \right\rangle \,,
\end{equation}
These are formally the same as the previously given expressions. The utility in this reformulation will be seen subsequently when the model is discretized. 

\subsection{Westervelt's equation and its Hamiltonian structure}

Westervelt's equation may be written
\begin{equation}
	p_{tt} - c^2 \Delta p - b \Delta p_t = (k p^2)_{tt}
\end{equation}
where $p$ is the pressure, $c$ is the speed of sound, $b$ controls the diffusion, and $k$ controls the strength of the nonlinearity. One begins by rewriting the dissipation-free part of the model as a first-order system in order to discern its Hamiltonian structure. For a dissipation-free model, let $b=0$. The first-order dissipation-free model is written as follows:
\begin{equation}
	\left[ (1 - k p )p \right]_{tt} - c^2 \Delta p = 0 \,.
\end{equation}
Letting $\bm{v}_t = - \nabla p$, this is easily rewritten as a first order system:
\begin{equation} \label{eq:westervelt_first_order}
	\begin{cases}
		[(1 - k p )p]_t + c^2 \nabla \cdot \bm{v} = 0 & \\
		\bm{v}_t + \nabla p = 0 \,. &
	\end{cases}
\end{equation}
This first order form then implies the appropriate constitutive relation: $(\bm{m}, \rho) = (\bm{v}, (1 - k p) p)$. One then infers that the auxiliary energy functional is given by
\begin{equation}
    K = \frac{k}{3} \int_\Omega p^3 \mathsf{d}^3 \bm{x} \,,
\end{equation}
and thus the Hamiltonian functional may be written 
\begin{multline}
    H[\rho, \bm{v}] 
    =
    K - \int_\Omega \frac{\delta K}{\delta p} p \mathsf{d}^3 \bm{x} + \frac{1}{2} \int_\Omega \left( p^2 + | \bm{v} |^2 \right) \mathsf{d}^3 \bm{x} \\
    = 
    \frac{1}{2} \int_\Omega \left[ \left( 1 - \frac{2 k p}{3} \right) p^2 + c^2 | \bm{v} |^2 \right] \mathsf{d}^3 \bm{x} \,,
\end{multline}
where $\rho = (1 - k p) p$. From the general theory, one knows that
\begin{equation}
    \frac{\delta H}{\delta \rho} = p \,,
    \quad \text{and} \quad
    \frac{\delta H}{\delta \bm{v}} = c^2 \bm{v} \,.
\end{equation}
Hence, using the general Poisson bracket for nonlinear acoustic wave models, \eqref{eq:poisson_bracket_general_acoustic}, one recovers the dissipation-free Westervelt equations. Dissipation is added to the model in the usual manner using the symmetric bracket given in equation \ref{eq:acoustic_symmetric_bracket}. It is straightforward to show that the evolution rule $\dot{F} = \{F, H\} + (F, H)$ recovers Westervelt's equations. 

If one rescales $x \to L \tilde{x}$, and $t \to T \tilde{t}$ such that $c = L/T$, then rescaling the parameters $b$ and $k$ appropriately, one obtains a dimensionless model:
\begin{equation}
    [(1 - k p) p]_{tt} - \Delta p - b \Delta p_t = 0 \,.
\end{equation}
This ensures that the wave speed is unity. This dimensionless model will be used for numerical experiments. For the purpose of convergence studies, it is helpful to add external forcing. This allows one to verify the method using the method of manufactured solutions. One adds forcing via
\begin{equation} \label{eq:westervelt_first_order_forcing}
	\begin{cases}
		[(1 - k p )p]_t + \nabla \cdot \bm{v} - b \Delta p = S_p & \\
		\bm{v}_t + \nabla p = \bm{S}_{\bm{v}} \,, &
	\end{cases}
\end{equation}
where $(S_p, \bm{S}_{\bm{v}})$ are right-hand source terms. This forcing term does not fit into the self-consistent Hamiltonian and dissipative structure and should be thought of as an external input to the system. 

\section{A discrete Hamiltonian theory of nonlinear acoustics} \label{sec:disc_ham_theory}

As a prelude to considering the mimetic discretization of Westervelt's equation, this section describes a general structure-preserving framework for nonlinear acoustic models, which encompasses finite element and finite difference approaches. This general approach suggests many other related discretization strategies for such models which might be studied in future works. 

\subsection{A general framework for mimetic discretization}

The differential operators appearing in nonlinear acoustic wave equations are the standard vector calculus operators which fit into the de Rham sequence \cite{frankel_2011, ArnoldDouglasN2010Feec}. It is helpful to briefly review structure-preserving grid-based discretizations of the de Rham sequence as these techniques will be used subsequently. Discrete representations of the de Rham sequence which preserve the de Rham cohomology take a general form which may be summarized by the following commuting diagram:
\begin{equation}
    \begin{tikzcd}[row sep=large]
    V^0 \simeq H^1(\Omega) 
    \arrow{r}{\msf{d}_0 = \bmsf{grad}} 
    \arrow{d}{\mathcal{R}^{0}} 
    \ar[dd, bend right=35, pos=0.2, swap, "\Pi^0"] & 
    V^{1} \simeq H(\bmsf{curl}, \Omega) 
    \arrow{r}{\msf{d}_1 = \bmsf{curl}} 
    \arrow{d}{\mathcal{R}^{1}} 
    \ar[dd, bend right=35, pos=0.2, swap, "\Pi^1"] &
    V^{2} \simeq H(\msf{div}, \Omega) 
    \arrow{r}{\msf{d}_2 = \msf{div}} 
    \arrow{d}{\mathcal{R}^{2}} 
    \ar[dd, bend right=35, pos=0.2, swap, "\Pi^2"] &
    V^{3} \simeq L^2(\Omega) 
    \arrow{d}{\mathcal{R}^{3}} 
    \ar[dd, bend right=35, pos=0.2, swap, "\Pi^3"] \\
    \mathcal{C}_{0} 
    \ar[r, crossing over, "\mathbbm{d}_{0} = \mathbb{G}"] 
    \arrow{d}{\mathcal{I}^{0}} & 
    \mathcal{C}_{1} 
    \ar[r, crossing over, "\mathbbm{d}_{1} = \mathbb{C}"] 
    \arrow{d}{\mathcal{I}^{1}} & 
    \mathcal{C}_{2} 
    \ar[r, crossing over, "\mathbbm{d}_{2} = \mathbb{D}"] 
    \arrow{d}{\mathcal{I}^{2}} &
    \mathcal{C}_{3} 
    \arrow{d}{\mathcal{I}^{3}} \\
    V^{0}_h 
    \arrow{r}{\msf{d}_{0,h} = \bmsf{grad}_h} & 
    V^{1}_h 
    \arrow{r}{\msf{d}_{1,h} = \bmsf{curl}_h} & 
    V^{2}_h 
    \arrow{r}{\msf{d}_{2,h} = \msf{div}_h} &
    V^{3}_h \,.
    \end{tikzcd} 
\end{equation}
The vector spaces making up this diagram are summarized as follows:
\begin{itemize}
    \item Top level: the continuous function spaces, $V^\ell$, making up the de Rham sequence in three-dimensions;
    \item Bottom level: the Galerkin subspaces, $V_h^\ell \subset V^\ell$;
    \item Middle level: coefficient spaces, $\mathcal{C}_\ell = \mathbb{R}^{N_\ell}$.
\end{itemize}
The operators in this diagram are as follows:
\begin{itemize}
    \item Reduction operators: $\mathcal{R}^\ell: V^\ell \to \mathcal{C}_\ell$;
    \item Interpolation operators: $\mathcal{I}^\ell: \mathcal{C}_\ell \to V^\ell_h$, defined as right inverses of the reduction operators, $\mathcal{R}^\ell \circ \mathcal{I}^\ell = \msf{id}$;
    \item Projection operators: $\Pi^\ell = \mathcal{I}^\ell \circ \mathcal{R}^\ell: V^\ell \to V^\ell_h$;
    \item Continuous exterior derivatives: $\msf{d}_\ell : V^\ell \to V^{\ell+1}$;
    \item Discrete exterior derivatives: $\msf{d}_{\ell,h} : V^\ell_h \to V^{\ell+1}_h$ such that
    \begin{equation}
        \msf{d}_{\ell,h} = \left. \msf{d}_\ell \right|_{V_h^\ell} \,;
    \end{equation}
    \item Derivative matrices: $\mathbbm{d}_\ell : \mathcal{C}_\ell \to \mathcal{C}_{\ell+1}$ such that
    \begin{equation}
    \begin{aligned}
        &\mathbbm{d}_\ell \mathcal{R}^\ell u = \mathcal{R}^{\ell+1} \msf{d}_\ell u \,, 
        \quad \forall u \in V^\ell \,, \\
        \quad \text{and} \quad 
        &\mathcal{I}^{\ell+1} \mathbbm{d}_\ell \bmsf{u} = \msf{d}_\ell \mathcal{I}^{\ell}  \bmsf{u} \,, 
        \quad \forall \bmsf{u} \in \mathcal{C}_\ell \,.
    \end{aligned}
    \end{equation}
\end{itemize}
Commutativity of the diagram ensures that the nilpotency of the exterior derivative is retained at the discrete level. Many discretization strategies, including mimetic finite differences \cite{bochev_and_hyman} and finite element exterior calculus \cite{ArnoldDouglasN2010Feec}, fit into this general paradigm with their particular properties being determined by the choice of reduction and interpolation operators. The particular mimetic finite difference method used in this work is described in appendix \ref{appendix:mimetic_fd}.

\subsection{Discrete duality pairings and discrete functional calculus} \label{sec:discrete_duality}

Let $V$ be a function space and let $V_h \subseteq V$ be a finite dimensional Galerkin subspace. Moreover, let $\mathcal{C} = \mathbb{R}^N$ be a space of Galerkin coefficients, and define the reduction and interpolation maps respectively as follows
\begin{equation}
    \mathcal{R}: V \to \mathcal{C} \,,
    \quad \text{and} \quad
    \mathcal{I}: \mathcal{C} \to V_h 
\end{equation}
such that reduction is a right inverse of interpolation, $\mathcal{R} \circ \mathcal{I} = \msf{id}$. Finally, define the projection operator such that
\begin{equation}
    \Pi: V \to V_h = \mathcal{I} \circ \mathcal{R} \,.
\end{equation}
The projection is assumed to be an approximation of the identity in the operator norm on $V$: $ \| \msf{id} - \Pi \| = O(h^p)$. 

\begin{theorem} \label{thm:disc_func_deriv}
Let $K: V \to \mathbb{R}$ be an arbitrary functional and let $\msf{K} := K \circ \mathcal{I}: \mathcal{C} \to \mathbb{R}$ represent the discrete analog of the functional $K$. Moreover, define $\bmsf{u} = \mathcal{R} u$ and $\bmsf{v} = \mathcal{R} v$. Then 
\begin{equation}
    \begin{aligned}
    \left\langle D(K \circ \Pi)[u], v \right\rangle_{V^*, V} 
    &= \left\langle D(K \circ \mathcal{I} )[ \mathcal{R} u ], \mathcal{R} v \right\rangle_{\mathcal{C}^*, \mathcal{C}} \\
    &= \left\langle D\msf{K}[ \bmsf{u} ], \bmsf{v} \right\rangle_{\mathcal{C}^*, \mathcal{C}} 
    =\left( \frac{\partial \msf{K}}{\partial \bmsf{u}} \right)^T \bmsf{v}.
    \end{aligned}
\end{equation}
Therefore, $D(K \circ \Pi): V \to \mathbb{R}$ is defined by the map
\begin{equation}
    D(K \circ \Pi) (\cdot) = \left( \frac{\partial \msf{K}}{\partial \bmsf{u}} \right)^T \mathcal{R} (\cdot) \,.
\end{equation}
\end{theorem}
This result points to a general representation of elements of the dual space $V_h^*$. Given a basis $\{ \Lambda_i \}_{i=1}^N \subseteq V_h$, its dual basis is made up of the individual reduction maps which return each degree of freedom:
\begin{equation}
    \mathcal{R} \circ \mathcal{I} = \msf{id} \implies
    \mathcal{R}_i \Lambda_j = \delta_{ij} \,.
\end{equation}
Define a coefficient space of dual degrees of freedom, $\mathcal{C}^* \simeq \mathcal{C}$. One may define a dual interpolation operator $\mathcal{I}^*: \mathcal{C}^* \to V^*_h$ such that
\begin{equation}
    \mathcal{I}^* \bmsf{u} = \bmsf{u}^T \mathcal{R}(\cdot) \,.
\end{equation}
Hence, using the evaluation pairing, one finds that $\langle \mathcal{I}^* \bmsf{u}, \mathcal{I} \bmsf{v} \rangle = \bmsf{u}^T \bmsf{v}$. One may think of the interpolation operator as a dual reduction operator, $\mathcal{R}^*: V^* \to \mathcal{C}^*$, defined such that $(\mathcal{R}^* v)_i = \langle v, \Lambda_i \rangle$. 
\begin{prop}
With $\mathcal{I}^*$ and $\mathcal{R}^*$ defined as above, it follows that
\begin{equation}
    \mathcal{R}^* \mathcal{I}^* = \langle \mathcal{R}(\cdot), \mathcal{I}(\cdot) \rangle = \msf{id} \in \mathcal{C}^* \times \mathcal{C}^* \,.
\end{equation}
\end{prop}
\noindent Finally, define the dual projection to be $\Pi^* = \mathcal{I}^* \circ \mathcal{R}^*$. This allows one to obtain a very similar result to theorem \ref{thm:disc_func_deriv}:
\begin{theorem} \label{thm:dual_disc_func_deriv}
Let $K: V^* \to \mathbb{R}$ be an arbitrary functional and let $\msf{K} := K \circ \mathcal{R}^*: \mathcal{C}^* \to \mathbb{R}$ represent the discrete analog of the functional $K$. Moreover, define $\bmsf{u} = \mathcal{R}^* u$ and $\bmsf{v} = \mathcal{R}^* v$. Then 
\begin{equation}
    \begin{aligned}
    \left\langle D(K \circ \Pi^*)[u], v \right\rangle_{V, V^*} 
    &= \left\langle D(K \circ \mathcal{I}^* )[ \mathcal{R}^* u ], \mathcal{R}^* v \right\rangle_{\mathcal{C}, \mathcal{C}^*} \\
    &= \left\langle D\msf{K}[ \bmsf{u} ], \bmsf{v} \right\rangle_{\mathcal{C}, \mathcal{C}^*} 
    =\left( \frac{\partial \msf{K}}{\partial \bmsf{u}} \right)^T \bmsf{v}.
    \end{aligned}
\end{equation}
Therefore, $D(K \circ \Pi^*): V^* \to \mathbb{R}$ is defined by the map
\begin{equation}
    D(K \circ \Pi^*) (\cdot) = \left( \frac{\partial \msf{K}}{\partial \bmsf{u}} \right)^T \mathcal{R}^* (\cdot) \,.
\end{equation}
But $D(K \circ \Pi^*)(\cdot) \in V^{**} \sim V$ assuming $V$ is reflexive (one will generally assume $V$ is Hilbert). In particular, one may identify
\begin{equation}
    D(K \circ \Pi^*)(\cdot) = \mathcal{I} \frac{\partial \msf{K}}{\partial \bmsf{u}} \,.
\end{equation}
\end{theorem}

This yields a key result which facilitates the discretization of Poisson brackets.
\begin{theorem}
For $F: V \to \mathbb{R}$ and $G: V^* \to \mathbb{R}$, one finds that
\begin{equation}
    \left\langle D(F \circ \Pi) (\cdot), D(G \circ \Pi^*)(\cdot) \right\rangle_{V^*, V}
    =
    D \msf{F}(\cdot)^T D \msf{G}(\cdot) \,,
\end{equation}
where $\msf{F} = F \circ \mathcal{R}: \mathcal{C} \to \mathbb{R}$ and $\msf{G} = G \circ \mathcal{R}^*: \mathcal{C}^* \to \mathbb{R}$.
\end{theorem}

Finally, a suitable discretization of the duality map $\mathcal{H}: V \to V^*$ must be prescribed. This is most reasonably given by
\begin{equation}
    \mathbb{H}: \mathcal{C} \to \mathcal{C}^* = \mathcal{R}^* \circ \mathcal{H} \circ \mathcal{I} \,.
\end{equation}
The precise form of this map depends on the manner in which the dual space $V^*$ is identified. In this work, lowest order finite differences (or equivalently, finite elements) are used. See appendix \ref{appendix:mimetic_fd} for further discussion of mimetic finite difference methods. 

\subsection{A mimetic discretization of nonlinear acoustics}

Suppose that one discretizes the fields via the Ansatz that $p_h \in V^0_h$, and $\bm{v}_h \in V^1_h$. That is, there exist $\bmsf{p} \in \mathcal{C}_0$ and $\bmsf{v} \in \mathcal{C}_1$ such that
\begin{equation}
    p_h = \mathcal{I}^0 \bmsf{p} \,,
    \quad \text{and} \quad
    \bm{v}_h = \mathcal{I}^1 \bmsf{v} \,.
\end{equation}
What distinguishes different nonlinear acoustic models is the auxiliary energy function, $K: V^0 \times V^1 \to \mathbb{R}$. A discrete analog is defined as follows:
\begin{equation}
    \msf{K}: \mathcal{C}_0 \times \mathcal{C}_1 \ni (\bmsf{p}, \bmsf{v}) 
    \mapsto
    K[\mathcal{I}^0 \bmsf{p}, \mathcal{I}^1 \bmsf{v}] \in \mathbb{R} \,.
\end{equation}
In some cases, it may be advantageous to further approximate this map by using an appropriate quadrature rule to compute $K[\mathcal{I}^0 \bmsf{p}, \mathcal{I}^1 \bmsf{v}]$. This consideration must tailored to each problem to find the most expeditious solution. The fields $\rho_h \in V_{0,h}^*$ and $\bm{m}_h \in V_{1,h}^*$ are defined such that
\begin{equation}
    \rho_h = \bm{\uprho}^T \mathcal{R}^0(\cdot) \,,
    \quad \text{and} \quad
    \bm{m}_h = \bmsf{m}^T \mathcal{R}^1(\cdot) \,.
\end{equation}
Hence, it follows that the discrete constitutive relations are given by
\begin{equation}
    \bm{\uprho} = \mathbb{H}_0 \bmsf{p} - \frac{\partial \msf{K}}{\partial \bmsf{p}} \,,
    \quad \text{and} \quad
    \bmsf{m} = \mathbb{H}_1 \bmsf{v} + \frac{\partial \msf{K}}{\partial \bmsf{v}} \,.
\end{equation}
The Hamiltonian is then given by
\begin{equation}
    \msf{H}[\bm{\uprho}, \bmsf{v}] = \msf{K}(\bmsf{p}, \bmsf{v}) 
    - \left( \frac{\partial \msf{K}}{\partial \bmsf{p}} \right)^T \bmsf{p} 
    + \frac{1}{2} 
    \left(
    \bmsf{p}^T \mathbb{H}_0 \bmsf{p}
    +
    \bmsf{v}^T \mathbb{H}_1 \bmsf{v}
    \right) \,,
\end{equation}
where $\bmsf{p} = \bmsf{p}(\bm{\uprho}, \bmsf{v})$ is defined implicitly through a coordinate change. One may show that
\begin{equation}
	\frac{\partial \msf{H}}{\partial \bm{\uprho}} = \bmsf{p}
	\quad \text{and} \quad
	\frac{\partial \msf{H}}{\partial \bmsf{v}} = \bmsf{m} \,,
\end{equation} 
where the same algebraic arguments may be used as in the continuous setting. 

The Poisson bracket is obtained by discretizing the functional derivatives according to theorems \ref{thm:disc_func_deriv} and \ref{thm:dual_disc_func_deriv}. One finds
\begin{equation}
    \{ \msf{F}, \msf{G} \}
    =
    \left( \frac{\partial \msf{F}}{\partial \bmsf{v}} \right)^T \mathbb{G}
    \frac{\partial \msf{G}}{\partial \bm{\uprho}}
    -
    \left( \frac{\partial \msf{G}}{\partial \bmsf{v}} \right)^T \mathbb{G}
    \frac{\partial \msf{F}}{\partial \bm{\uprho}} \,,
\end{equation}
where $\mathbb{G}$ is the discrete gradient matrix implied by the mimetic discretization framework. Similarly, the discrete symmetric bracket is given by
\begin{equation}
    ( \msf{F}, \msf{G} )
    =
    - b
    \left( \mathbb{G} \frac{\partial \msf{F}}{\partial \bm{\uprho}} \right)^T 
    \mathbb{H}_1
    \mathbb{G}
    \frac{\partial \msf{G}}{\partial \bm{\uprho}} \,.
\end{equation}
The equations of motion are obtained by noting that an arbitrary observable $\msf{F}$ evolves as $\dot{\msf{F}} = \{\msf{F}, \msf{H}\} + (\msf{F}, \msf{H})$. This implies
\begin{equation} \label{eq:disc_nl_acoustic}
    \frac{\msf{d}}{\msf{d} t}
    \begin{pmatrix}
        \bm{\uprho} \\
        \bmsf{v}
    \end{pmatrix}
    =
    \begin{pmatrix}
        - b \mathbb{G}^T \mathbb{H}_1 \mathbb{G} & \mathbb{G}^T \\
        - \mathbb{G} & 0 
    \end{pmatrix}
    \begin{pmatrix}
        \bmsf{p} \\
        \bmsf{m}
    \end{pmatrix}
    \,,
\end{equation}
where
\begin{equation} \label{eq:disc_cons_law}
    \bm{\uprho} = \mathbb{H}_0 \bmsf{p} - \frac{\partial \msf{K}}{\partial \bmsf{p}} \,,
    \quad \text{and} \quad
    \bmsf{m} = \mathbb{H}_1 \bmsf{v} + \frac{\partial \msf{K}}{\partial \bmsf{v}} \,.
\end{equation}
This system has a discrete energy dissipation rate
\begin{equation}
    \dot{\msf{H}} = - b \bmsf{p}^T \mathbb{G}^T \mathbb{H}_1 \mathbb{G} \bmsf{p} \,,
\end{equation}
and possesses the conservation law $\mathbb{C} \dot{\bmsf{v}} = 0$. 

\section{A Hamiltonian structure-preserving discretization of Westervelt's equation} \label{sec:disc_westervelt}

The general spatial discretization framework for nonlinear acoustic models described in the previous section is now applied to Westervelt's equation.

\subsection{A mimetic discretization of Westervelt's equation}

For the sake of simplicity, the mimetic discretization framework uses lowest order interpolation on uniform, equispaced grid on a periodic domain, $\Omega = [0,1]^d$. This ensures that the duality matrices are diagonal and the derivative matrices are circulant, see appendix \ref{appendix:mimetic_fd}. On a periodic domain, it is advantageous to use trapezoidal rule to approximate integrals due to the exponential convergence property. That is, one defines $\msf{K}: \mathcal{C}^0 \to \mathbb{R}$ such that
\begin{equation}
    \msf{K}(\bmsf{p}) = \frac{k}{3} \sum_{i = 1}^{N_0} \msf{p}_i^3 h^d \,,
\end{equation}
where $d$ is the number of spatial dimensions, and $h$ is the grid spacing. The validity of this approximation follows because Lagrange interpolation (used to interpolate $0$-forms) is interpolatory at the grid points. Any other structure-preserving collocation method would likewise allow the Hamiltonian to be approximated via the trapezoidal rule. Contrast this with a structure-preserving finite element method which would yield a discrete energy functional of the form
\begin{equation}
    \msf{K}(\bmsf{p}) 
    =
    \frac{k}{3} \sum_{i,j,k = 1}^{N_0} \msf{p}_i \msf{p}_j \msf{p}_k \int_\Omega \Lambda^0_i(\bm{x}) \Lambda^0_j(\bm{x}) \Lambda^0_k(\bm{x}) \mathsf{d}^3 \bm{x} \,.
\end{equation}
In general, this coupling tensor may be quite complicated. While it is possible to efficiently compute the entries of this rank-$3$ tensor via quadrature methods, the nonlinear constitutive relation implied by equation \eqref{eq:disc_cons_law} generically yields a global coupling of all degrees of freedom. Hence, a collocation method is preferable for algorithmic expedience. 

The discrete Hamiltonian implied by a collocation method with the nonlinear functional, $K$, approximated via trapezoidal rule is given by
\begin{equation}
    \msf{H}(\bmsf{p}, \bmsf{v})
    =
    \frac{1}{2} \left( 
    \bmsf{p}^T \mathbb{H}_0 \bmsf{p}
    +
    \bmsf{v}^T \mathbb{H}_1 \bmsf{v}
    \right) - \frac{2 k}{3} \sum_{i = 1}^{N_0} \msf{p}_i^3 h^d \,.
\end{equation}
It is well-known that trapezoidal rule achieves exponential convergence for smooth fields on a periodic domain \cite{doi:10.1137/130932132}: hence, this is an effective choice of quadrature rule. The constitutive relations take the form
\begin{equation}
    \bm{\uprho} = \mathbb{H}_0 \bmsf{p} - k \bmsf{p} \pdot \bmsf{p} h^d \,,
    \quad \text{and} \quad
    \bmsf{m} = \mathbb{H}_1 \bmsf{v} \,,
\end{equation}
where $\pdot$ denotes the pointwise (Hadamard) product. If $\mathbb{H}_0$ is diagonal, i.e.\ if an orthogonal basis is used, then the nonlinear solve decouples to yield exactly solvable scalar nonlinear constitutive relations for $\msf{p}_i = \msf{p}_i(\uprho_i)$. This work only considers one- and two-dimensional problems with lowest-order finite mimetic difference methods, a particular case in which the Hodge duality matrices are diagonal. See appendix \ref{appendix:mimetic_fd} for a discussion of these duality matrices and an explanation for the connection between these matrices and the Hodge star operator \cite{frankel_2011} which motivates the terminology. 

\subsection{Structure-preserving time integration}

The final consideration in designing a structure-preserving numerical method is discerning an appropriate time-stepping method. This is done with operator splitting techniques. In particular, individual partial flows are obtained for the dissipative and conservative dynamics before composing these together using Strang splitting \cite{doi:10.1137/0705041} to get a scheme for the full dynamics. 

\subsubsection{Partial flow for the conservative dynamics}

The conservative dynamics, obtained by setting $b=0$, are likewise integrated via a splitting method. In this case, the Hamiltonian is split into two pieces:
\begin{equation}
    \msf{H}_\rho =
    \msf{K}(\bmsf{p}) 
    - \left( \frac{\partial \msf{K}}{\partial \bmsf{p}} \right)^T \bmsf{p} 
    + \frac{1}{2} 
    \bmsf{p}^T \mathbb{H}_0 \bmsf{p} \,,
    \quad \text{and} \quad
    \msf{H}_v = \frac{1}{2} \bmsf{v}^T \mathbb{H}_1 \bmsf{v} \,.
\end{equation}
In the case of lowest order interpolation, one has
\begin{equation}
    \msf{H}_\rho = \sum_{i = 1}^{N_0} \left( 1 - \frac{2 k}{3} \msf{p}_i^2 \right) \msf{p}_i h^d \,,
\end{equation}
which causes the nonlinear constitutive relation to completely decouple into independent scalar problems. This splits the dynamics into two decoupled exactly solvable subsystems. These yield independent update rules for $\bmsf{p}$ and $\bmsf{v}$. The rule for $\bmsf{p}$ is as follows:
\begin{equation}
    \varphi_{\rho}(\Delta t) : (\bm{\uprho}, \bmsf{v}) 
        \mapsto (\bm{\uprho} + \Delta t \mathbb{G}^T \mathbb{H}_1 \bmsf{v}, \bmsf{v}) \,.
\end{equation}
The rule for $\bmsf{v}$ is as follows:
\begin{equation}
    \varphi_{v}(\Delta t) : (\bm{\uprho}, \bmsf{v}) \mapsto (\bm{\uprho}, \bmsf{v} - \Delta t \mathbb{G} \bmsf{p}(\bm{\uprho})) \,,
\end{equation}
where $\bmsf{p}(\bm{\uprho})$ solves
\begin{equation}
    \bm{\uprho} = \mathbb{H}_0 \bmsf{p} - \frac{\partial \msf{K}}{\partial \bmsf{p}} \,,
    \quad \text{or} \quad
    \bm{\uprho} = ( \bmsf{p} - k \bmsf{p} \pdot \bmsf{p} ) h^d 
\end{equation}
in the case of lowest order interpolation. It is helpful to maintain a state vector which stores both $\bm{\uprho}$ and $\bmsf{p}$. One then must update $\bmsf{p}$ whenever $\bm{\uprho}$ is updated. Hence, in practice, a nonlinear solve is performed after every $\varphi_{\rho}$ step. 

An order $O(\Delta t)$ flow is obtained by composing together the two flows using Lie-Trotter splitting \cite{trotter1959product}:
\begin{equation}
    \varphi_{cons}^{(1)}(\Delta t) = \varphi_{v}(\Delta t) \circ \varphi_{\rho}(\Delta t) \,.
\end{equation}
An order $O(\Delta t^2)$ flow is obtained by composing together the two flows symmetrically using Strang splitting \cite{doi:10.1137/0705041}:
\begin{equation}
    \varphi_{cons}^{(2)}(\Delta t) = \varphi_{v}(\Delta t/2) \circ \varphi_{\rho}(\Delta t) \circ \varphi_{v}(\Delta t/2) \,.
\end{equation}
The order does not impact the order of convergence, although it is better to take only a single step with $\varphi_{\rho}$ since this partial step involves a nonlinear solve. 

If one wishes to include external forcing, it is natural to bundle this part of the evolution with the conservative time-step for the density. This is because this time-step may be accomplished exactly. Letting
\begin{equation}
    \bmsf{S}_p(t) = \mathcal{R}^0(S_p) \,,
    \quad \text{and} \quad
    \bmsf{S}_v(t) = \mathcal{R}^1(\bm{S}_{\bm{v}}) \,,
\end{equation}
the contribution of the external forcing to the $n^{th}$ time-step may be computed exactly:
\begin{equation}
\begin{aligned}
    &\varphi_{\rho}(t_n, \Delta t) : 
        (\bm{\uprho}, \bmsf{v}) 
        \mapsto 
        \left(
        \bm{\uprho} 
        + \Delta t \mathbb{G}^T \mathbb{H}_1 \bmsf{v} 
        + \int_{t_n}^{t_{n+1}} \mathbb{H}_0 \bmsf{S}_p(t) \mathsf{d} t, 
        \bmsf{v}
        \right) 
        \,, \\
    \quad \text{and} \quad
    &\varphi_{v}(t_n, \Delta t) : 
        (\bm{\uprho}, \bmsf{v}) 
        \mapsto 
        \left(
        \bm{\uprho}, \bmsf{v} 
        - \Delta t \mathbb{G} \bmsf{p}(\bm{\uprho})
        + \int_{t_n}^{t_{n+1}} \bmsf{S}_v(t) \mathsf{d} t
        \right)
        \,.
\end{aligned}
\end{equation}
It is necessary to multiply the pressure source term by $\mathbb{H}_0$ because this equation is modeled on the dual space. It is also worth noting that it is necessary to know the current time to apply the time-advance map with time-dependent external forcing. Finally, assuming the integral of the right-hand side forcing is simple enough that its integral may be computed exactly, this partial flow may be computed exactly. 

\subsubsection{Partial flow for the dissipative dynamics} 

The dissipative dynamics involve a heat-like relaxation on $\bm{\uprho}$. This otherwise simple evolution is complicated by the nonlinearity of the model. It is not possible to prescribe an exact flow for the dissipative dynamics as was done in the conservative case. Rather, the dissipative flow is approximated via a Runge-Kutta method. Because the full flow is approximated with second-order Strang splitting, a second-order method is needed

It is desirable to design this relaxation such that the correct rate of energy dissipation is achieved. That is, one wishes for energy to dissipate as
\begin{equation}
    \left. \dot{H} \right|_{\bmsf{p} = \bmsf{p}^n}
    =
    \left. \left(\frac{\partial \msf{H}}{\partial \bm{\uprho}} \right)^T \dot{\bm{\uprho}} \right|_{\bmsf{p} = \bmsf{p}^n}
    =
    - b (\bmsf{p}^n)^T \mathbb{G}^T \mathbb{H}_1 \mathbb{G} \bmsf{p}^n
     \,,
\end{equation}
or an approximation thereof. Here, time derivatives describe evolution under the dissipative dynamics only. 

An $O(\Delta t)$ approximate flow is obtained via forward Euler:
\begin{equation}
    \bm{\uprho}^{n+1} = \bm{\uprho}^{n} - b \Delta t \mathbb{G}^T \mathbb{H}_1 \mathbb{G} \bmsf{p}^n \,. 
\end{equation}
Hence, the update rule is
\begin{equation}
    \varphi_{diss}^{(1)}(\Delta t) : (\bm{\uprho}, \bmsf{v}) 
        \mapsto (\bm{\uprho} - b \Delta t \mathbb{G}^T \mathbb{H}_1 \mathbb{G} \bmsf{p}(\bm{\uprho}), \bmsf{v}) \,.
\end{equation}
An $O(\Delta t^2)$ approximate flow is obtained via the midpoint method (i.e.\ second-order Runge-Kutta):
\begin{equation}
    \bm{\uprho}^{n+1} = \bm{\uprho}^{n} - b \Delta t \mathbb{G}^T \mathbb{H}_1 \mathbb{G} 
    \bmsf{p}\left( 
    \bm{\uprho}^{n} - \frac{b \Delta t}{2} \mathbb{G}^T \mathbb{H}_1 \mathbb{G} \bmsf{p}(\bm{\uprho}^n)
    \right) 
    \,. 
\end{equation}
Hence, the update rule is
\begin{equation}
    \varphi_{diss}^{(2)}(\Delta t) : (\bm{\uprho}, \bmsf{v}) 
    \mapsto 
    \left(
    \bm{\uprho} 
    - b \Delta t \mathbb{G}^T \mathbb{H}_1 \mathbb{G} 
    \bmsf{p} \left( 
    \bm{\uprho} - \frac{b \Delta t}{2} \mathbb{G}^T \mathbb{H}_1 \mathbb{G} \bmsf{p}(\bm{\uprho})
    \right) , \bmsf{v}
    \right) \,.
\end{equation}
As in the conservative dynamics, one must solve a nonlinear system to get $\bmsf{p}(\bm{\uprho})$. The second-order method requires that one solve the nonlinear system twice. Hence, it is critical that an efficient nonlinear solver be used. 

\subsubsection{Full time-step}

An $O(\Delta t)$ flow may be obtained by composing together the two flows:
\begin{equation}
    \varphi_{LT}(\Delta t) = \varphi_{diss}^{(1)}(\Delta t) \circ \varphi_{cons}^{(1)}(\Delta t) \,.
\end{equation}
As each stage is $O(\Delta t)$, the overall method is $O(\Delta t)$. A second-order method may be accomplished with Strang splitting \cite{doi:10.1137/0705041}:
\begin{multline}
    \varphi_S(\Delta t) 
    = 
    \varphi_{cons}^{(2)}(\Delta t/2)^* 
    \circ \varphi_{diss}^{(2)}(\Delta t) 
    \circ \varphi_{cons}^{(2)}(\Delta t/2) \\
    =
    \varphi_{v}(\Delta t/4)
    \circ \varphi_{\rho}(\Delta t/2) 
    \circ \varphi_{v}(\Delta t/4) 
    \circ \varphi_{diss}^{(2)}(\Delta t) \\
    \circ \varphi_{v}(\Delta t/4) 
    \circ \varphi_{\rho}(\Delta t/2) 
    \circ \varphi_{v}(\Delta t/4) \,.
\end{multline}
Each stage in the splitting method must be a second-order method for the overall method to be second-order. After steps involving $\varphi_\rho$, one must solve for $\bmsf{p}$, and each dissipative step requires two nonlinear solves. For the second-order method, a total of four nonlinear solves are needed per time-step. 

\subsubsection{Nonlinear solve for the pressure}

In this work, the spatial discretization is done using a collocation method because of the nonlinear system relating $\bmsf{p}$ with $\bm{\uprho}$ decouples into independent scalar problems. One wishes to invert
\begin{equation}
    \bm{\uprho}(\bmsf{p}) = h^d (1 - k \bmsf{p}) \pdot \bmsf{p} = \bm{\uprho}^\star \,,
\end{equation}
or, component-wise,
\begin{equation}
    \rho_i^\star = h^d(1 - k p_i) p_i \,.
\end{equation}
The inverse is explicitly computable:
\begin{equation}
    p_i^{\pm}
    =
    \frac{1 \pm \sqrt{1 - \frac{4 k \rho_i^\star}{h^d}}}{2 k} \,.
\end{equation}
One must determine which branch is appropriate. One can see that, as $k \to 0$, 
\begin{equation}
    p_i^+ \sim \frac{1}{k} - \frac{\rho_i^\star}{h} + O(k) \,,
    \quad \text{while} \quad
    p_i^- \sim \frac{\rho_i^\star}{h} + O(k) \,.
\end{equation}
Hence, the negative branch corresponds with physical solutions. Thus, the update rule is given by
\begin{equation}
    \bmsf{p} = \frac{1 - \sqrt{1 - \frac{4 k \bm{\uprho}^\star}{h^d}}}{2 k} \,,
\end{equation}
where the algebra is performed component-wise. Because the nonlinear solve may be performed analytically when linear interpolation is used, the cost of the nonlinear solve is insignificant. 

\subsubsection{Conditions for stability of the time-step}

Each time-step is composed of conservative and dissipative partial flows. Each of these pieces of the full time-step individually needs to be stable. To begin, note that
\begin{equation}
    p(\rho) = \frac{1 - \sqrt{1 - 4 k \rho}}{2 k} 
    \implies
    p'(\rho) = \frac{1}{\sqrt{1 - 4 k \rho}} \,.
\end{equation}
For sufficiently small $k$ and $ \| \rho_0 \|_{\infty} \ll 1/(4k)$, the stability of the time-advance map of the nonlinear system is guaranteed by ensuring that the linear problem obtained when $k \to 0$ is stable. Hence, it is sufficient to ensure the stability of the linearized time-advance map, see appendix \ref{appendix:nonlinear_stability} for further discussion. 

For the conservative part of the dynamics, one needs that 
\begin{equation}
    \Delta t \sqrt{ \| c_S^2 \|_{\infty} \| \mathbb{H}_0^{-1} \mathbb{G}^T \mathbb{H}_1 \mathbb{G} \| } < 1 \,,
\end{equation}
where $c_S^2(\bm{x})$ is the sound speed, and for the dissipative dynamics, one requires that
\begin{equation}
    \Delta t \| b \|_{\infty} \| \mathbb{H}_0^{-1} \mathbb{G}^T \mathbb{H}_1 \mathbb{G} \| < 1 \,.
\end{equation}
For generality, the sound speed and dissipation parameter are allowed to be functions of space. Hence, one chooses
\begin{equation}
    \Delta t < \min \left\{ \frac{1}{ \sqrt{ \| c_S^2 \|_{\infty} \| \mathbb{H}_0^{-1} \mathbb{G}^T \mathbb{H}_1 \mathbb{G} \| }}, \frac{1 }{\| b \|_{\infty} \| \mathbb{H}_0^{-1} \mathbb{G}^T \mathbb{H}_1 \mathbb{G} \| } \right\} \,.
\end{equation}
Due to the simplicity of the finite difference stencils at lowest order, the stability is easily ascertained from the Gershgorin circle theorem \cite{izvestija/gerschgorin31}. At lowest-order in one dimension, the finite difference stencil is given by
\begin{equation}
    \mathbb{H}_0^{-1} \mathbb{G}^T \mathbb{H}_1 \mathbb{G}
    =
    \frac{1}{\Delta x^2} 
    \begin{pmatrix}
        -1 & 2 & -1
    \end{pmatrix}
    \implies
    \| \mathbb{H}_0^{-1} \mathbb{G}^T \mathbb{H}_1 \mathbb{G} \|
    < \frac{4}{\Delta x^2} \,.
\end{equation}
The higher-dimensional Laplacian matrices are simply Kronecker sums of one-dimensional Laplacian matrices. Hence, in two-dimensions
\begin{equation}
    \| \mathbb{H}_0^{-1} \mathbb{G}^T \mathbb{H}_1 \mathbb{G} \|
    <
    4 \left( \frac{1}{\Delta x^2} + \frac{1}{\Delta y^2} \right) \,.
\end{equation}
One can see that, generally, the time-step size is limited by the stability condition for the dissipative dynamics. The method might be improved by solving the dissipative sub-step implicitly. Any implicit method for the dissipative dynamics would require solving a nonlinear system which couples adjacent degrees of freedom for the pressure and is not considered in this work. 

\section{Numerical examples} \label{sec:numerical_examples}

In the following numerical examples, two particular qualities of this finite-difference solver for Westervelt's equation are of particular interest. The first quality of interest is the convergence of the method: this is investigated using the method of manufactured solutions. The second quality of interest is its structure-preserving properties, e.g.\ achieving the proper energy dissipation rate and conserving vorticity. Tests are performed in one and two spatial dimensions. In all examples, second order methods are used for temporal integration. 

\subsection{One-dimensional examples}

In one spatial dimension, the lowest order duality matrices are given by
\begin{equation}
    \mathbbm{h}_0 = \Delta x \mathbb{I} \,,
    \quad \text{and} \quad
    \mathbbm{h}_1 = \Delta x^{-1} \mathbb{I} \,,
\end{equation}
where $h$ is the grid spacing and $\mathbb{I}$ is the identity matrix. The duality matrices are written in lower-case blackboard font for notational convenience as they will be used to construct the duality matrices in two-dimensions as well. The derivative matrix is circulant:
\begin{equation}
    \mathbbm{d}_0 
    =
    \begin{pmatrix}
        -1 & 1 & 0 & \hdots & 0 \\
        0 & -1 & 1 & 0 & \\
        \vdots & & \ddots & \ddots & \\
        0 & & & -1 & 1 \\
        1 & 0 & \hdots & 0 & -1
    \end{pmatrix} \,.
\end{equation}
These choices are consistent with a lowest order mimetic finite difference method. Notice, this choice of duality and derivative matrix results in a derivative operator which is identical to a standard lowest-order finite difference stencil. Given a grid resolution, $\Delta x$, the time-step is taken to be $\Delta t = \Delta x^2/ (4 b)$. 

\subsubsection{Convergence study}

A convergence study is performed using the method of manufactured solutions. If one lets
\begin{equation}
    p(x,t) = \cos( 2 \pi t) \cos(2 \pi x) \,,
    \quad \text{and} \quad
    v(x,t) = \sin(2 \pi t) \sin(2 \pi x) \,,
\end{equation}
then it follows that
\begin{equation}
    S_p(x,t) = 4 \pi [ b \pi + k \sin(2 \pi t) \cos(2 \pi x) ] \cos(2 \pi t) \cos(2 \pi x) \,,
\end{equation}
while $S_v(x,t) = 0$. Moreover, one should use the initial data
\begin{equation}
    p_0(x) = \cos(2 \pi x) \,,
    \quad \text{and} \quad
    v_0(x) = 0 \,.
\end{equation}
The solution is computed on the spatio-temporal domain $\Omega \times [0,T] = [0,1] \times [0,1]$ with parameters
\begin{equation}
    k = 2 \times 10^{-1} \,,
    \quad \text{and} \quad
    b = 10^{-2} \,.
\end{equation}
The $L^2(\Omega \times [0,T])$ and $L^\infty(\Omega \times [0,T])$ norms of the error between the computed and exact solution given in figure \ref{table:west1d}. The errors are given relative to the exact solution:
\begin{equation}
    \msf{rel}_p(p_{num}, p_{exact})
    =
    \frac{\| p_{num} - p_{exact} \|_{L^p(\Omega \times [0,T])}}{\| p_{exact} \|_{L^p(\Omega \times [0,T])}} \,.
\end{equation}
One can see that second-order convergence is achieved as is predicted by the theory. 

\begin{table}
\centering
\begin{tabular}{||c c c c c||} 
 \hline
 $N_x$ & $L^2$ error & order & $L^\infty$ error & order \\ [0.5ex] 
 \hline\hline
 20 & 0.2385 & --- & 0.3281 & --- \\ 
 \hline
 40 & 0.05774 & 2.046 & 0.07817 & 2.069 \\
 \hline
 80 & 0.01433 & 2.011 & 0.01930 & 2.018 \\
 \hline
 160 & 0.003574 & 2.003 & 0.004808 & 2.005 \\
 \hline
 320 & 0.0008932 & 2.000 & 0.001201 & 2.001 \\ [1ex] 
 \hline
\end{tabular}
\caption{The relative error in the manufactured solution and the estimated order of convergence is given.}
\label{table:west1d}
\end{table}

\subsubsection{Gaussian wave propagation}

\begin{figure}
    \centering
    \includegraphics[width=\linewidth]{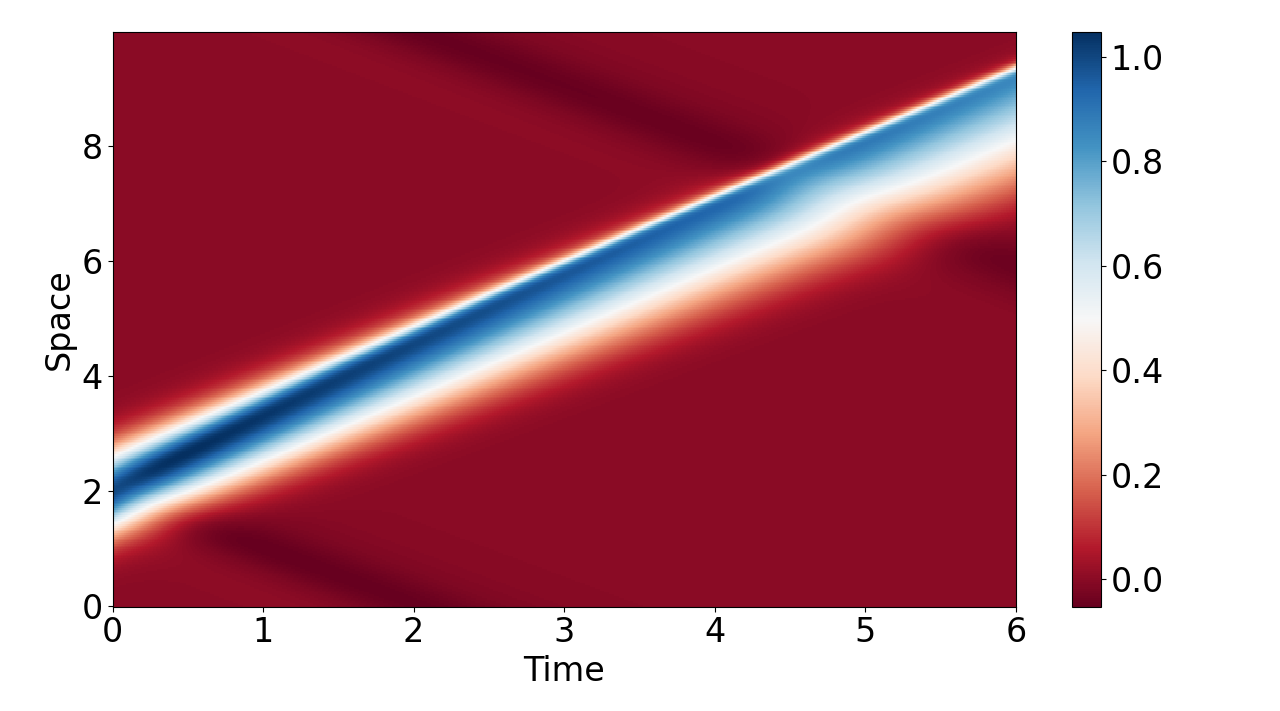}
    \caption{The evolution of a Gaussian wave-form predicted by the one-dimensional solver.}
    \label{fig:gaussian_wave_1d}
\end{figure}

A useful qualitative test case in one spatial dimension is the propagation of a Gaussian waveform:
\begin{equation}
    p_0(x) = \exp \left( - \frac{1}{2} \frac{(x/L - \mu)^2}{\sigma^2} \right) \,,
    \quad \text{and} \quad
    v_0(x) = 0 \,,
\end{equation}
where $L$ is the characteristic length of the domain. The spatio-temporal domain is taken to be $\Omega \times [0,T] = [0,10] \times [0,6]$. The parameters for the initial data are
\begin{equation}
    \mu = 0.2 \,,
    \quad
    L = 10 \,,
    \quad \text{and} \quad
    \sigma = 0.05 \,.
\end{equation}
The physical parameters are taken to be
\begin{equation}
    k = 2 \times 10^{-1} \,,
    \quad \text{and} \quad
    b = 10^{-2} \,.
\end{equation}
The domain is discretized with $N_x = 320$ grid-points. See figure \ref{fig:gaussian_wave_1d} for a visualization of the solution. One can see that the wave to steepens as it propagates, a common feature of nonlinear wave models.  

\subsubsection{The energy dissipation rate}

\begin{figure}
    \centering
    \includegraphics[width=\linewidth]{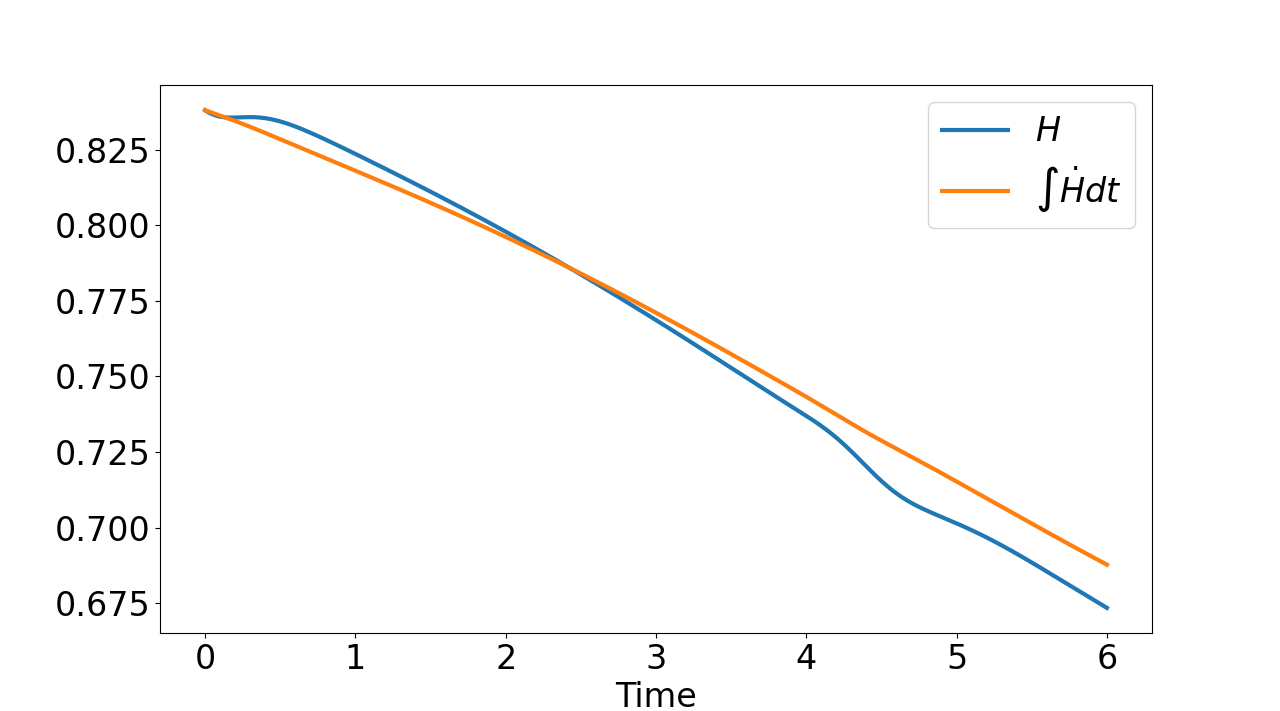}
    \caption{A comparison of the actual and predicted energy dissipation rate from the one-dimensional solver.}
    \label{fig:energy_dissipation_1d}
\end{figure}

Energy is dissipated at a known rate in the spatially-discrete model:
\begin{equation}
    \left. \dot{\msf{H}} \right|_{\bmsf{p} = \bmsf{p}^n}
    =
    - b (\bmsf{p}^n)^T \mathbb{G}^T \mathbb{H}_1 \mathbb{G} \bmsf{p}^n
     \,.
\end{equation}
In the following examples, this property is verified numerically. See figure \ref{fig:energy_dissipation_1d} for a comparison of the computed energy as a function of time, and the integral of the estimated energy dissipation rate as a function of time from the Gaussian wave propagation example. The two curves roughly coincide, but the non-conservation of energy due to the symplectic integrator for the conservative dynamics results in noticeable differences. If one were to use an energy conserving method for the conservative dynamics, it is likely that the two curves would better coincide. 

\subsection{Two-dimensional examples}

In two spatial dimensions, on a tensor product domain all discrete operators are built as tensor products of the one-dimensional operators. See \cite{10.1007/s10915-022-01781-3} for a description of the construction in the context of finite element methods on tensor product domains. The finite difference method used here is essentially just the lowest order finite element method. The duality matrices are given by
\begin{equation}
    \mathbbm{H}_0 = \mathbbm{h}_0^x \otimes \mathbbm{h}_0^y \,,
    \quad \text{and} \quad
    \mathbbm{H}_1 = 
    \begin{pmatrix}
        \mathbbm{h}_1^x \otimes \mathbbm{h}_0^y & 0 \\
        0 & \mathbbm{h}_0^x \otimes \mathbbm{h}_1^y
    \end{pmatrix} \,,
\end{equation}
where $\mathbbm{h}_0^x$ and $\mathbbm{h}_0^y$ are the $0$-form duality matrices for the one-dimensional $x$- and $y$-grids respectively. These are diagonal at lowest order. The gradient matrix is given by
\begin{equation}
    \mathbb{G}
    =
    \begin{pmatrix}
        \mathbbm{d}_0 \otimes \mathbb{I} & 0 \\
        0 & \mathbb{I} \otimes \mathbbm{d}_0 
    \end{pmatrix} \,.
\end{equation}
The discrete two-dimensional curl operator approximates the differential operator
\begin{equation}
    \msf{curl} \bm{v}
    =
    \partial_y v_x - \partial_x v_y \,,
\end{equation} 
and is given by
\begin{equation}
    \mathbb{C} =
    \begin{pmatrix}
        \mathbb{I} \otimes \mathbbm{d}_0 & 
        - \mathbbm{d}_0 \otimes \mathbb{I}
    \end{pmatrix} \,.
\end{equation}
Notice, $\mathbb{C} \mathbb{G} \equiv 0$. Given a grid resolution in each direction, $(\Delta x, \Delta y)$, the time-step is taken to be 
\begin{equation}
    \Delta t = \left( 4 \left( \frac{1}{ \Delta x^2} + \frac{1}{\Delta y^2} \right) b \right)^{-1} \,.
\end{equation}

\subsubsection{Convergence study}
A convergence study is performed using the method of manufactured solutions. If one lets
\begin{equation}
    p(x,y,t) 
    =
    \cos(2 \pi t) \cos(2 \pi x) \cos(2 \pi y) \,,
\end{equation}
and
\begin{equation}
    \bm{v}(x,y,t)
    =
    \begin{pmatrix}
        \sin(2 \pi t) \sin(2 \pi x) \cos(2 \pi y) \\
        \sin(2 \pi t) \cos(2 \pi x) \sin(2 \pi y) 
    \end{pmatrix} \,,
\end{equation}
then
\begin{multline}
    S_p(x,y,t)
    =
    2 \pi 
    \left[  
    4 b \pi \cos(2 \pi t) 
    + \sin(2 \pi t)
    + k \cos(2 \pi x) \cos(2 \pi y) \sin(4 \pi t) 
    \right] \\
    \times \cos(2 \pi x) \cos(2 \pi y) 
    \,,
\end{multline}
while $\bm{S}_{\bm{v}} = 0$. The solution is computed on the spatio-temporal domain $\Omega \times [0,T] = [0,1]^2 \times [0,1]$ with parameters
\begin{equation}
    k = 2 \times 10^{-1} \,,
    \quad \text{and} \quad
    b = 10^{-2} \,.
\end{equation}
The relative error in the $L^2(\Omega \times [0,T])$ and $L^\infty(\Omega \times [0,T])$ norms is given in figure \ref{table:west2d}. One can see that approximately second-order convergence is achieved as is predicted by the theory.

\begin{table}
\centering
\begin{tabular}{||c c c c c||} 
 \hline
 $N_x \times N_y$ & $L^2$ error & order & $L^\infty$ error & order \\ [0.5ex] 
 \hline\hline
 $20 \times 20$ & 0.1650 & --- & 0.2376 & --- \\ 
 \hline
 $40 \times 40$ & 0.04277 & 1.948 & 0.06575 & 1.853 \\
 \hline
 $80 \times 80$ & 0.01078 & 1.988 & 0.01689 & 1.961 \\
 \hline
 $160 \times 160$ & 0.002702 & 1.996 & 0.004250 & 1.991 \\
 \hline
 $320 \times 320$ & 0.0006758 & 1.999 & 0.001064 & 1.998 \\ [1ex] 
 \hline
\end{tabular}
\caption{The relative error in the manufactured solution and the estimated order of convergence is given.}
\label{table:west2d}
\end{table}

\subsubsection{Wave propagation in nonuniform medium} \label{sec:nonunif_medium}

A medium with spatially varying sound speed may be described with a simple modification of the Hamiltonian:
\begin{equation}
    H[\rho, \bm{v}]
    =
    \frac{1}{2} \int_\Omega \left[ \left( 1 - \frac{2 k p}{3} \right) p^2 + c^2_S(\bm{x}) | \bm{v} |^2 \right] \mathsf{d}^3 \bm{x} \,,
\end{equation}
where, as usual, $p = p(\rho)$ implicitly by inverting $\rho(p) = (1 - kp)p$. This variable sound speed is a consequence of variable background density in the acoustic medium. At lowest order interpolation, a non-uniform sound speed might be accommodated by letting
\begin{equation}
    \msf{H}_v = \frac{1}{2} \bmsf{v}^T \mathbbm{c}_S^2 \bmsf{v} \,,
\end{equation}
where one defines the matrix
\begin{equation}
    \mathbbm{c}_S^2
    =
    \tilde{\mathcal{R}}_1
    (
    \star
    c_S^2
    \mathcal{I}_1
    )
    \,.
\end{equation}
where $\tilde{\mathcal{R}}_1$ is the $1$-form reduction operator on the dual grid, see appendix \ref{appendix:mimetic_fd}. At lowest order, this matrix is diagonal with
\begin{equation}
    \mathbbm{c}_S^2
    =
    \begin{pmatrix}
        \mathbbm{c}_{S,x}^2 & 0 \\
        0 & \mathbbm{c}_{S,y}^2
    \end{pmatrix}
\end{equation}
where, recalling that $i = (i_x, i_y)$ is a multi-index,
\begin{equation}
    (\mathbbm{c}_{S,x}^2)_{ii}
    =
    \frac{1}{\Delta x} \int_{x_{i_x-1/2}}^{x_{i_x+1/2}} c^2_S(x, y_{i_y+1/2}) \mathsf{d} x \,,
    \quad \text{and} \quad
    (\mathbbm{c}_{S,y}^2)_{ii}
    =
    \frac{1}{\Delta y} \int_{y_{i_y-1/2}}^{y_{i_y+1/2}} c^2_S(x_{i_x+1/2}, y) \mathsf{d} x \,,
\end{equation}
where $x_{i_x \pm 1/2} = x_{i_x} \pm \Delta x/2$ and $y_{i_y \pm 1/2} = y_{i_y} \pm \Delta y/2$. That is, the kinetic energy is weighted by the local values of the sound speed. To leading order, this is well approximated by
\begin{equation}
    (\mathbbm{c}_{S,x}^2)_{ii}
    \approx
    c^2_S(x_{i_x}, y_{i_y+1/2}) \,,
    \quad \text{and} \quad
    (\mathbbm{c}_{S,y}^2)_{ii}
    \approx
    c^2_S(x_{i_x+1/2}, y_{i_y}) \,.
\end{equation}
This further approximation will be used in the numerical examples in this section. 

\begin{figure}
    \centering
    \includegraphics[width=0.9\linewidth]{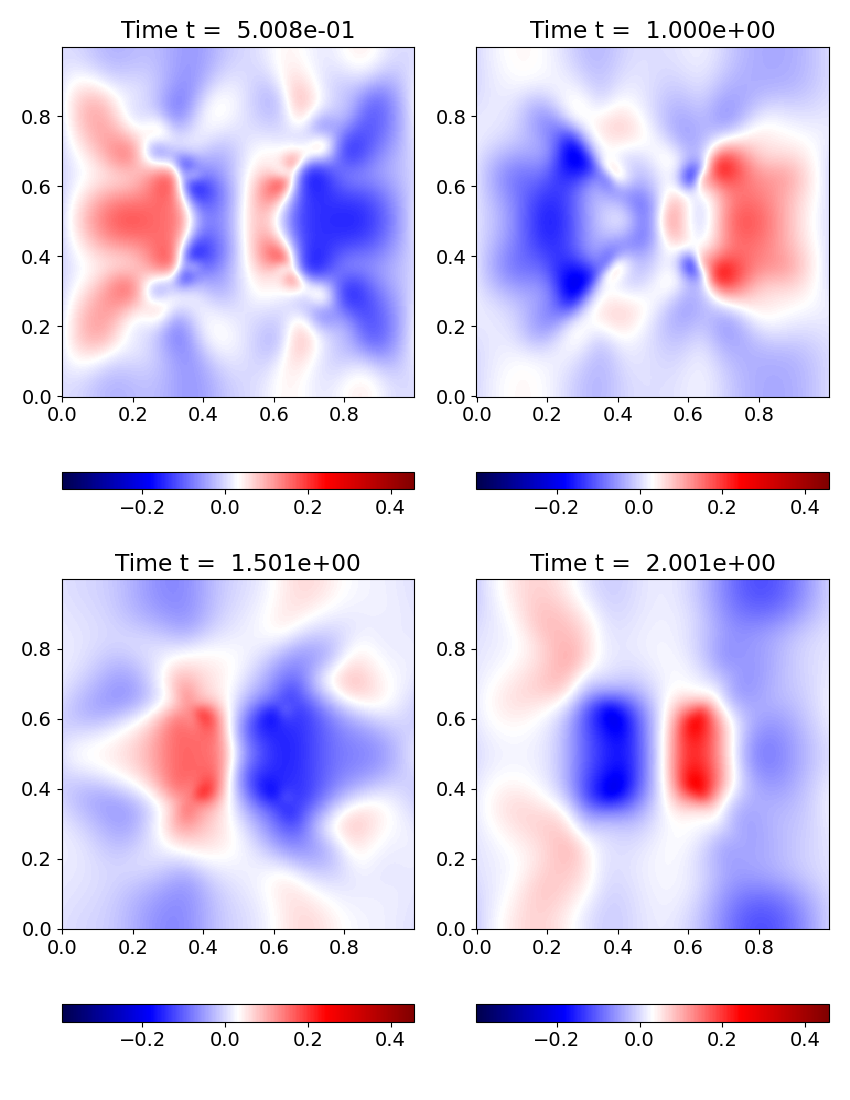}
    \caption{Visualization of the solution to Westervelt's equation in two dimensions with spatially varying sound speed.}
    \label{fig:west2d_snapshots}
\end{figure}

In the following test, the sound speed is taken to be
\begin{equation}
    c_S^2(x,y)
    =
    1 
    - \frac{3}{4} \sin^2(\pi x) \sin^2(\pi y)
    - \frac{1}{2} \sin^2(2 \pi x) \sin^2(2 \pi y)
    - \frac{1}{4} \sin^2(4 \pi x) \sin^2(4 \pi y) \,.
\end{equation}
This choice of sound-speed takes values in
\begin{equation}
    c_S^2(\bm{x}) \in [0.0782421, 1.0] \,,
    \quad \text{for} \quad
    \bm{x} \in [0,1]^2 \,.
\end{equation}
This function is chosen to allow for significant variation in the sound-speed in the domain. The initial conditions are chosen to be
\begin{equation}
    p_0(x,y) = 0 \,,
    \quad
    \bm{v}(x,y)
    =
    \begin{pmatrix}
        \cos(2 \pi x) \cos(2 \pi y) \\
        \sin(2 \pi x) \sin(2 \pi y) 
    \end{pmatrix} \,.
\end{equation}
These, relatively non-physical, initial conditions are chosen so that the vorticity is nonzero:
\begin{equation}
    \omega = \msf{curl} \bm{v}
    = 
    \partial_x v_y - \partial_y v_x 
    =
    4 \pi \cos(2 \pi x) \sin(2 \pi y)
    \neq 0 \,.
\end{equation}
This is chosen to verify the conservation of vorticity throughout the simulation. The physical parameters are
\begin{equation}
    k = 2 \times 10^{-1}, 
    \quad \text{and} \quad
    b = 10^{-3} \,.
\end{equation}
The spatial resolution is chosen to be $N_x \times N_y = 320 \times 320$, and the temporal resolution is $\Delta t = (4 b (N_x^2 + N_y^2))^{-1}$. See figure \ref{fig:west2d_snapshots} for a visualization of the solution. 

\subsubsection{The energy dissipation rate and conservation of vorticity}

\begin{figure}
    \centering
    \begin{subfigure}
        \centering
        \includegraphics[width=0.86\textwidth]{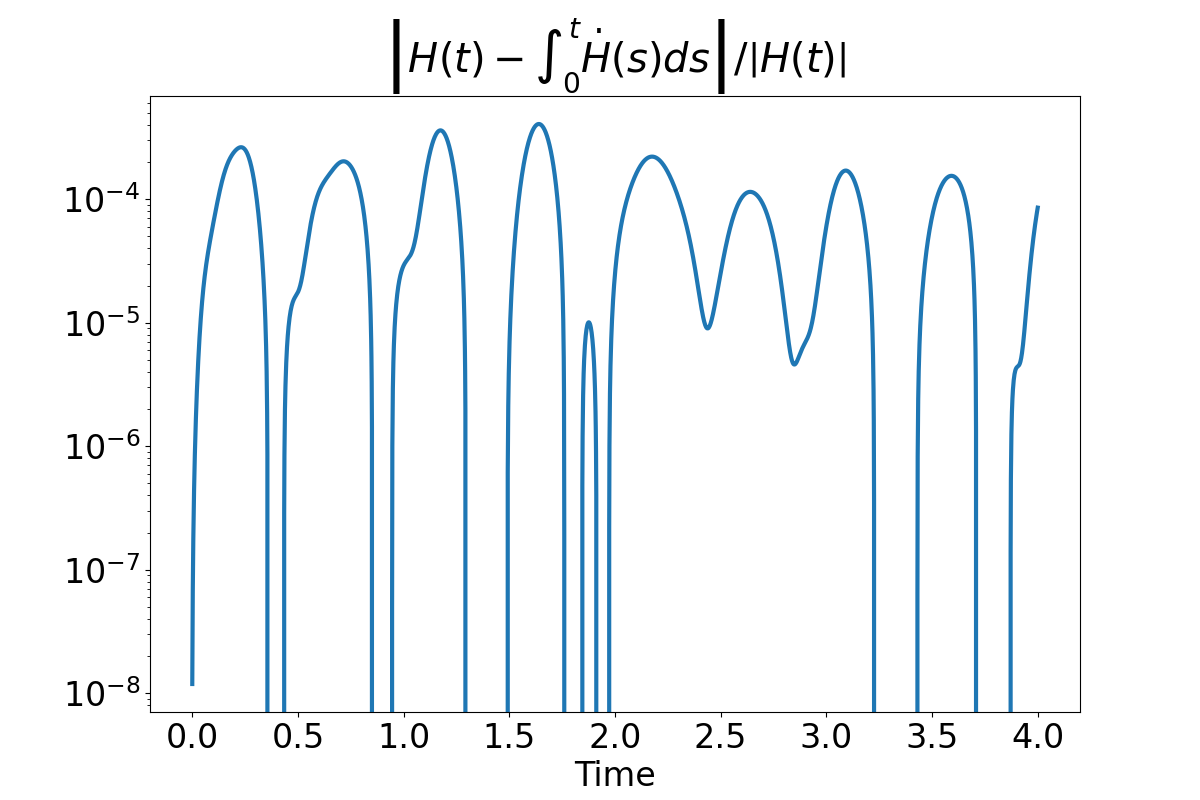}
    \end{subfigure}%
    \begin{subfigure}
        \centering
        \includegraphics[width=0.86\textwidth]{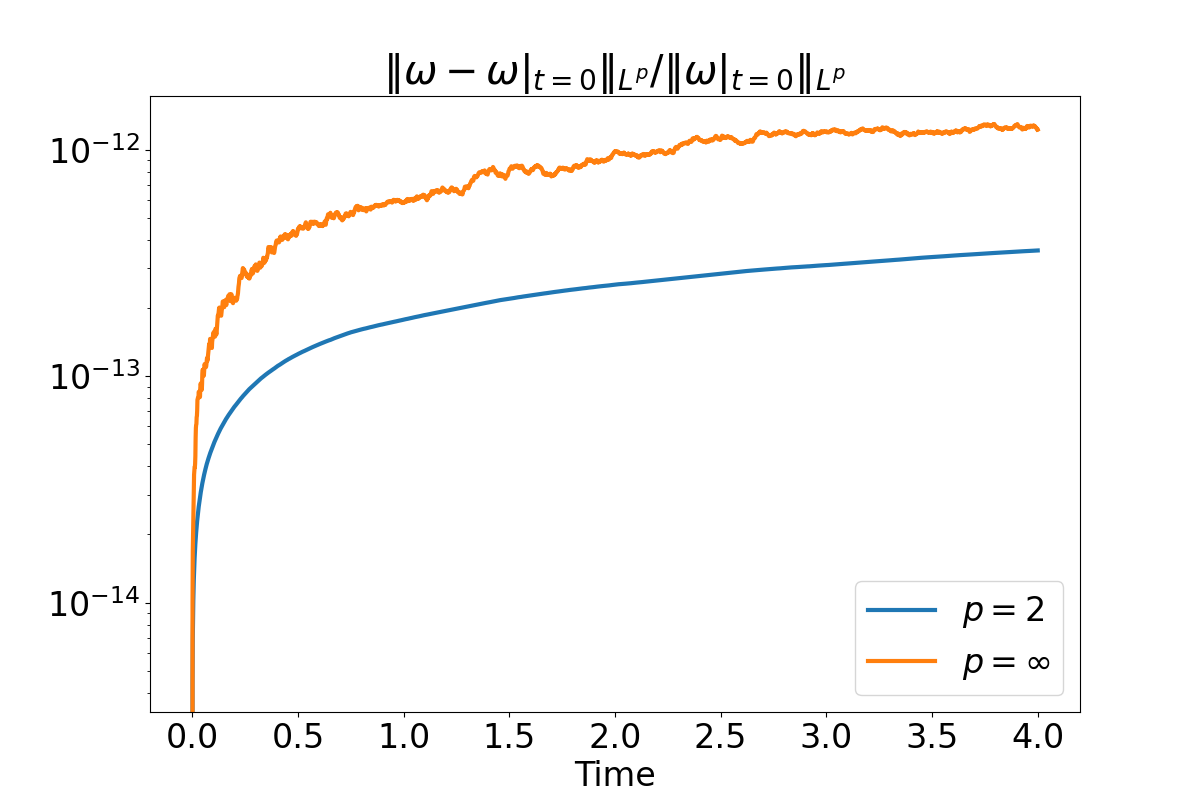}
    \end{subfigure}
    \caption{Scalar diagnostics for two-dimensional solver. Top: relative error between the actual and predicted energy dissipation rate. Bottom: relative error in the vorticity field as a function of time.}
    \label{fig:west2d_diags}
\end{figure}

The vorticity, that is $\bm{\omega} = \nabla \times \bm{v}$, is conserved as a Casimir invariant of the Hamiltonian structure. This may be easily verified by noting that 
\begin{equation}
    \partial_t \bm{v} + \nabla p = 0 
    \implies
    \partial_t \nabla \times \bm{v} = 0 \,.
\end{equation}
The spatially discrete model likewise conserves vorticity since, in a structure-preserving spatial discretization, $\mathbb{C} \mathbb{G} = 0$:
\begin{equation}
    \partial_t \bmsf{v} + \mathbb{G} \bmsf{p} = 0
    \implies
    \partial_t \mathbb{C} \bmsf{v} = 0 \,.
\end{equation}
As in the one-dimensional case, energy is dissipated at a known rate in the spatially-discrete model:
\begin{equation}
    \left. \dot{\msf{H}} \right|_{\bmsf{p} = \bmsf{p}^n}
    =
    - b (\bmsf{p}^n)^T \mathbb{G}^T \mathbb{H}_1 \mathbb{G} \bmsf{p}^n
     \,.
\end{equation}
See figure \ref{fig:west2d_diags} for a visualization of the rate of energy dissipation, and the conservation of vorticity from the numerical example in section \ref{sec:nonunif_medium}. As expected, the vorticity is conserved approximately to machine precision, while the energy is dissipated at approximately the correct rate modulo oscillation from the symplectic integrator. 

\section{Conclusion} \label{sec:conclusion}

This work presented a novel structure-preserving discretization strategy for Westervelt's equation based on its Hamiltonian structure. The spatial discretization yields a Hamiltonian ordinary differential equation in its dissipation-free limit and possesses a discrete energy-dissipation rate in the finite-dissipation regime. Moreover, the Casimir invariants of the continuous model are preserved by the discrete theory. The temporal discretization is accomplished via splitting methods. The method applies quite generally to many nonlinear acoustic models, but only Westervelt's equation was considered in detail due to its relative simplicity. Moreover, the spatial discretization accommodates a broad class of structure-preserving discretization approaches provided they preserve the de Rham cohomology. While a mimetic finite difference approach was used, one could also use a FEEC method \cite{ArnoldDouglasN2010Feec}, an appropriate structure-preserving spectral method \cite{10.1007/s10915-022-01781-3}, or structure-preserving pseudo-spectral methods \cite{robidoux_2008}. Even the simple, lowest-order finite difference approximation used herein was found to yield an efficient, convergent method which well approximates the true energy dissipation rate and exactly conserves the curl of the velocity field. This work only presented the method as a proof of concept leaving substantial room for future inquiry. Future work will consider other nonlinear acoustic models, Galerkin and spectral methods, higher-order methods, more general domains, and physically meaningful problems. 

\bibliographystyle{plain} 		 
\bibliography{references}

\appendix

\section{Mimetic finite differences as a special case of mimetic discretization} \label{appendix:mimetic_fd}

The mimetic finite difference approximation used in this work was only briefly sketched in the body of the work. Here, further detail is provided. See \cite{bochev_and_hyman} and references therein for additional information on the mimetic finite difference method. 

\subsection{The discrete double de Rham complex}

One can explicitly construct a discrete double de Rham complex based on staggered grids \cite{PALHA20141394, kreeft2011mimetic}. Let $\tilde{V}^\ell$ denote the space of twisted differential forms \cite{eldred_and_bauer, frankel_2011}. The de Rham complex of straight and twisted differential forms are related via the Hodge star operator. In three dimensions, one writes
\begin{equation}
\begin{tikzcd}
    V^0 \ar[r, "\msf{d}_{0}"] \ar[d, shift left, "\star_0"] 
    & V^1 \ar[r, "\msf{d}_{1}"] \ar[d, shift left, "\star_1"]
    & V^2 \ar[r, "\msf{d}_{2}"] \ar[d, shift left, "\star_2"]
    & V^3 \ar[d, shift left, "\star_3"] \\
    \tilde{V}^3 \ar[u, shift left, "\tilde{\star}_3"]
    & \ar[l, "\msf{d}_{2}"] \tilde{V}^2 \ar[u, shift left, "\tilde{\star}_2"]
    & \ar[l, "\msf{d}_{1}"] \tilde{V}^1 \ar[u, shift left, "\tilde{\star}_1"]
    & \ar[l, "\msf{d}_{0}"] \tilde{V}^0 \ar[u, shift left, "\tilde{\star}_0"] \,.
\end{tikzcd}
\end{equation}
A discrete de Rham complex is defined analogously for both the straight and twisted forms, but the twisted forms are defined using a staggered dual grid. The notion of using a staggered dual grid to represent the dual de Rham complex has previously been explored in the literature \cite{PALHA20141394, kreeft2011mimetic}. The use of staggered grids likewise has a precedent in the finite difference literature \cite{doi:10.1137/0727052}. The use of staggered grids is straightforward with periodic boundary conditions as in this work. 

To be explicit, suppose $\Omega = [0,1]^3$. The discrete geometry for the primal grid is specified as a tensor-product grid with a uniform grid in each direction: 
\begin{equation}
    \{ (x_{i_1}, y_{i_2}, z_{i_3}) \}_{\bm{i} = 0}^{\bm{N}_0} \,,
    \quad \text{where} \quad 
    \bm{i} \in \llbracket 1, \bm{N}_0 \rrbracket \coloneq \prod_{\alpha = 1}^3 \llbracket 1, N_0^\alpha \rrbracket \,,
\end{equation}
where $\llbracket 1, N \rrbracket = \{1, 2, \hdots, N-1, N\}$, and one defines 
\begin{equation}
    x_i = \frac{i}{N_0^1} \in [0, 1]
\end{equation}
and analogously for $y$ and $z$. The dual grid is given by the midpoints of the primal grid in each direction: 
\begin{equation}
    \{ (x_{i_1 + 1/2}, y_{i_2 + 1/2}, z_{i_3 + 1/2}) \}_{\bm{i} = 0}^{\bm{N}_0} \,,
    \quad \text{where} \quad 
    \bm{i} \in \llbracket 1, \bm{N}_0 \rrbracket \coloneq \prod_{\alpha = 1}^3 \llbracket 1, N_0^\alpha \rrbracket \,,
\end{equation}
and the dual grid points are defined to be
\begin{equation}
    x_{i+1/2} = \frac{x_i + x_{i+1}}{2} \,,
\end{equation}
and similarly for $y$ and $z$. 

Because the three-dimensional discrete geometry is given by a tensor product of one-dimensional grids, the reduction operators for the three-dimensional de Rham complex are given by tensor products of the one-dimensional reduction operators. Letting $(x_i, y_j, z_k) \in \mathtt{N}$, the primal reduction operators are defined to be
\begin{itemize}
    \item $\mathcal{R}_0: V^0 \to \mathcal{C}^0$ associates a scalar-valued function with with its values at the vertices of the primal grid:
    \begin{equation}
        \mathcal{R}_{0,(i,j,k)}(\phi) = \phi(x_i, y_j, z_k) = \bm{\phi}_{i,j,k} \,.
    \end{equation}
    \item $\mathcal{R}_1: V^1 \to \mathcal{C}^1$ associates a vector-valued function with its circulations along the edges of the primal grid: 
    \begin{equation}
        \mathcal{R}_{1,(i,j,k)}(\bm{E}) = 
        \begin{pmatrix}
            \mathcal{R}^x_{1,(i,j,k)}(E_x) \\
            \mathcal{R}^y_{1,(i,j,k)}(E_y) \\
            \mathcal{R}^z_{1,(i,j,k)}(E_z)
        \end{pmatrix} \,,
    \end{equation}
    where
    \begin{equation}
        \mathcal{R}^x_{1,(i,j,k)}(E_x)
        \eqcolon
        \bmsf{E}_{i+1/2,j,k}
        \eqcolon
        \bmsf{E}^x_{i,j,k}
        \eqcolon
        \int_{x_i}^{x_{i+1}} E_x(x,y_j,z_k) \mathsf{d} x \,,
    \end{equation}
    and similarly for the remaining components. 
    \item $\mathcal{R}_2: V^2 \to \mathcal{C}^2$ associates a vector-valued function with its fluxes through the faces of the primal grid: 
    \begin{equation}
        \mathcal{R}_{2,(i,j,k)}(\bm{B}) = 
        \begin{pmatrix}
            \mathcal{R}^x_{2,(i,j,k)}(B_x) \\
            \mathcal{R}^y_{2,(i,j,k)}(B_y) \\
            \mathcal{R}^z_{2,(i,j,k)}(B_z)
        \end{pmatrix} \,,
    \end{equation}
    where
    \begin{equation}
        \mathcal{R}^x_{2,(i,j,k)}(B_x)
        \eqcolon
        \bmsf{B}_{i,j+1/2,k+1/2}
        \eqcolon
        \bmsf{B}^x_{i,j,k}
        \eqcolon
        \int_{y_j}^{y_{j+1}} \int_{z_k}^{z_{k+1}} B_x(x_i,y,z) \mathsf{d} y \mathsf{d} z \,,
    \end{equation}
    and similarly for the remaining components. 
    \item $\mathcal{R}_3: V^3 \to \mathcal{C}^3$ associates a scalar-valued function with with its cell integrals on the primal grid:
    \begin{equation}
        \mathcal{R}_{3,(i,j,k)}(\rho) 
        = \int_{x_i}^{x_{i+1}} \int_{y_j}^{y_{j+1}} \int_{z_k}^{z_{k+1}} \rho(x,y,z) \mathsf{d} x \mathsf{d} y \mathsf{d} z 
        = \bm{\rho}_{i,j,k} \,.
    \end{equation}
\end{itemize}
The reduction operators on the dual complex are defined analogously using the dual grid. 

Interpolation is accomplished via interpolation/histopolation \cite{10.1007/978-3-642-15337-2_17}. The basis functions are built from tensor products of one-dimensional basis functions with the property
\begin{equation}
    \ell_i(x_j) = \delta_{ij} \,,
    \quad \text{and} \quad
    \int_{x_j}^{x_{j+1}} \mathfrak{e}_i(x) \mathsf{d} x = \delta_{ij} \,.
\end{equation}
Evidently, $\{ \ell_i \}$ are the Lagrange interpolating polynomials while $\{ \mathfrak{e}_i \}$ are known as edge polynomials \cite{10.1007/978-3-642-15337-2_17}. It is possible to use high-order interpolating polynomials, but this work only considers linear interpolation. With these one may define the three-dimensional interpolation operators: 
\begin{itemize}
    \item $\mathcal{I}^0: \mathcal{C}^0 \to V^0_h$ is defined
    \begin{equation}
        \mathcal{I}^0 \bm{\upphi} 
        =
        \sum_{i,j,k}
        \upphi_{i,j,k} \ell_i(x) \ell_j(y) \ell_k(z)
        \,.
    \end{equation}
    \item $\mathcal{I}^1: \mathcal{C}^1 \to V^1_h$ is defined
    \begin{equation}
        \mathcal{I}^1 \bmsf{E}
        =
        \sum_{i,j,k}
        \begin{pmatrix}
            \msf{E}_{i,j,k}^x \mathfrak{e}_i(x) \ell_j(y) \ell_k(z) \\
            \msf{E}_{i,j,k}^y \ell_i(x) \mathfrak{e}_j(y) \ell_k(z) \\
            \msf{E}_{i,j,k}^z \ell_i(x) \ell_j(y) \mathfrak{e}_k(z) 
        \end{pmatrix}
        \,.
    \end{equation}
    \item $\mathcal{I}^2: \mathcal{C}^2 \to V^2_h$ is defined
    \begin{equation}
        \mathcal{I}^2 \bmsf{B}
        =
        \sum_{i,j,k}
        \begin{pmatrix}
            \msf{B}_{i,j,k}^x \ell_i(x) \mathfrak{e}_j(y) \mathfrak{e}_k(z) \\
            \msf{B}_{i,j,k}^y \mathfrak{e}_i(x) \ell_j(y) \mathfrak{e}_k(z) \\
            \msf{B}_{i,j,k}^z \mathfrak{e}_i(x) \mathfrak{e}_j(y) \ell_k(z) 
        \end{pmatrix}
        \,.
    \end{equation}
    \item $\mathcal{I}^3: \mathcal{C}^3 \to V^3_h$ is defined
    \begin{equation}
        \mathcal{I}^3 \bm{\uprho} 
        =
        \sum_{i,j,k}
        \uprho_{i,j,k} \mathfrak{e}_i(x) \mathfrak{e}_j(y) \mathfrak{e}_k(z)
        \,.
    \end{equation}
\end{itemize}
With these definitions, one can show that $\mathcal{R}^\ell \mathcal{I}^\ell = \msf{id}$. The interpolation operators on the dual complex are defined analogously over the dual grid. 

The gradient, curl, and divergence matrices follow directly from the definition of the reduction operators and the stipulation that
\begin{equation}
    \begin{aligned}
        \mathcal{R}_1( \bmsf{grad} \phi ) &= \mathbb{G} \mathcal{R}_0(\phi) \,,
            \quad \forall \phi \in V^0 \,, \\
        \mathcal{R}_2( \bmsf{curl} \bm{E} ) &= \mathbb{C} \mathcal{R}_1(\bm{E}) \,,
            \quad \forall \bm{E} \in V^1 \,, \\
        \mathcal{R}_3(\msf{div} \bm{B} ) &= \mathbb{D} \mathcal{R}_2(\bm{B}) \,,
            \quad \forall \bm{B} \in V^2 \,, \\
        \tilde{\mathcal{R}}_1( - \bmsf{grad} \phi ) &= \tilde{\mathbb{G}} \tilde{\mathcal{R}}_0(\phi) \,,
            \quad \forall \phi \in \tilde{V}^0 \,, \\
        \tilde{\mathcal{R}}_2( \bmsf{curl} \bm{E} ) &= \tilde{\mathbb{C}} \tilde{\mathcal{R}}_1(\bm{E}) \,,
            \quad \forall \bm{E} \in \tilde{V}^1 \,, \\
        \text{and} \quad
        \tilde{\mathcal{R}}_3(- \msf{div} \bm{B} ) &= \tilde{\mathbb{D}} \tilde{\mathcal{R}}_2(\bm{B}) \,,
            \quad \forall \bm{B} \in \tilde{V}^2 \,.
    \end{aligned}
\end{equation}
The one dimensional derivative matrix is defined to be the matrix such that
\begin{equation}
    \int_{x_i}^{x_{i+1}} \frac{\mathsf{d}}{\mathsf{d} x} \phi(x) \mathsf{d} x 
    = \phi(x_{i+1}) - \phi(x_i) 
    = (\mathbbm{d}_N \bm{\upphi})_i \,,
\end{equation}
where $\upphi_i = \phi(x_i)$. Due to the assumption of periodic boundary conditions, one finds that
\begin{equation}
    \mathbbm{d}_N
    =
    \begin{pmatrix}
        -1 & 1 & 0 & \hdots & 0 \\
        0 & -1 & 1 & 0 & \\
        \vdots & & \ddots & \ddots & \\
        0 & & & -1 & 1 \\
        1 & 0 & \hdots & 0 & -1
    \end{pmatrix}
    \in \mathbb{R}^{N \times N} \,.
\end{equation}
One may show that
\begin{itemize}
    \item the primal complex gradient matrix is given by
    \begin{equation}
        \mathbb{G} = 
        \begin{pmatrix}
            \mathbbm{d}_{N_1} \otimes \mathbb{I}_{N_2} \otimes \mathbb{I}_{N_3} \\
            \mathbb{I}_{N_1} \otimes \mathbbm{d}_{N_2} \otimes \mathbb{I}_{N_3} \\
            \mathbb{I}_{N_1} \otimes \mathbb{I}_{N_2} \otimes \mathbbm{d}_{N_3} 
        \end{pmatrix} \,.
    \end{equation}
    \item the primal complex curl matrix is given by
    \begin{equation}
        \mathbb{C}
        =
        \begin{pmatrix}
            \mathbb{O}_{N} 
                & - \mathbb{I}_{N_1} \otimes \mathbb{I}_{N_2} \otimes \mathbbm{d}_{N_3} 
                    & \mathbb{I}_{N_1} \otimes \mathbbm{d}_{N_2} \otimes \mathbb{I}_{N_3} \\
            \mathbb{I}_{N_1} \otimes \mathbb{I}_{N_2} \otimes \mathbbm{d}_{N_3} 
                & \mathbb{O}_{N} 
                    & - \mathbbm{d}_{N_1} \otimes \mathbb{I}_{N_2} \otimes \mathbb{I}_{N_3} \\
            - \mathbb{I}_{N_1} \otimes \mathbbm{d}_{N_2} \otimes \mathbb{I}_{N_3} 
                & \mathbbm{d}_{N_1} \otimes \mathbb{I}_{N_2} \otimes \mathbb{I}_{N_3}
                    & \mathbb{O}_{N}
        \end{pmatrix} \,,
    \end{equation}
    where $N = N_1 + N_2 + N_3$. The discrete curl matrix is self-adjoint: $\mathbb{C} = \mathbb{C}^T$ since $\mathbbm{d}_N^T = - \mathbbm{d}_N$.
    \item and the primal complex divergence matrix is given by
    \begin{equation}
        \mathbb{D}
        =
        \begin{pmatrix}
            \mathbbm{d}_{N_1} \otimes \mathbb{I}_{N_2} \otimes \mathbb{I}_{N_3} &
            \mathbb{I}_{N_1} \otimes \mathbbm{d}_{N_2} \otimes \mathbb{I}_{N_3} &
            \mathbb{I}_{N_1} \otimes \mathbb{I}_{N_2} \otimes \mathbbm{d}_{N_3} 
        \end{pmatrix}
        = - \mathbb{G}^T \,.
    \end{equation}  
\end{itemize}
The dual complex derivative matrices are derived via similar means. One can show that
\begin{equation}
    \mathbb{G}^T = \tilde{\mathbb{D}} \,,
    \quad
    \mathbb{C}^T = \tilde{\mathbb{C}} \,,
    \quad \text{and} \quad
    \mathbb{D}^T = \tilde{\mathbb{G}} \,.
\end{equation}
Other than specifying the discrete duality pairings, the reduction, interpolation, and discrete derivative operators together specify the entire discrete double de Rham complex. The construction in two-dimensions is analogous. 

\subsection{Duality and the discrete Hodge operator}

The choice of duality map, $\mathcal{H}: V^\ell \to (V^\ell)^*$, is a key feature in the design of a mimetic discretization. The problem of defining an effective discrete Hodge star operator has generated significant interest in the literature \cite{Hiptmair2001DiscreteHO, bochev_and_hyman, kreeft2011mimetic}. Two options are of particular note. If one uses $L^2$-duality to identify the dual space, the discrete duality operator, $\mathbb{H}_\ell: \mathcal{C}^\ell \to (\mathcal{C}^\ell)^*$, is the the Gram matrix formed from the basis functions used for interpolation. Using the definition of the interpolation and dual reduction operators, one finds that
\begin{equation}
    \bmsf{u}^T \mathbb{H}_\ell \bmsf{v}
    =
    \bmsf{u}^T \left( 
    (\mathcal{R}^\ell)^* \circ \mathcal{H}^\ell \circ \mathcal{I}^\ell
    \right)
    \bmsf{v}
    =
    \langle 
    \mathcal{H}^\ell \mathcal{I}^\ell \bmsf{v}, 
    \mathcal{I}^\ell \bmsf{u} 
    \rangle
    =
    (
    \mathcal{I}^\ell \bmsf{v}, 
    \mathcal{I}^\ell \bmsf{u}
    )_{L^2}
    =
    (
    \mathcal{I}^\ell \bmsf{u}, 
    \mathcal{I}^\ell \bmsf{v}
    )_{L^2}
    \,.
\end{equation}
Depending on the choice of interpolation operators, this choice coincides with a finite element method. 

Alternatively, one may specify duality through the double de Rham complex. If one makes the identification $(\mathcal{C}^\ell)^* \simeq \tilde{\mathcal{C}}^{n-\ell}$, then one may define discrete duality via
\begin{equation}
    \mathcal{H}^\ell \circ \mathcal{I}^\ell
    =
    \tilde{\mathcal{R}}^{n-\ell} \circ \star_\ell \,.
\end{equation}
This identification is possible if $\text{dim}(\mathcal{C}^\ell) = \text{dim}(\tilde{\mathcal{C}}^{n-\ell})$. With this identification yields a discrete duality map of the form
\begin{equation}
    \mathbb{H}_\ell 
    =
    (\mathcal{R}^\ell)^*
    \circ
    \tilde{\mathcal{R}}^{n-\ell} 
    \circ 
    \star_\ell
    =
    \tilde{\mathcal{R}}^{n-\ell} 
    \circ 
    \star_\ell
    \circ
    \mathcal{I}^\ell : \mathcal{C}^\ell \to \tilde{\mathcal{C}}^{n-\ell} \,.
\end{equation}
This is frequently called the natural Hodge discrete star operator \cite{bochev_and_hyman}. This approach to duality is frequently used in mimetic finite difference methods. One has that
\begin{equation}
    \tilde{\mathbb{H}}_{n-\ell} \mathbb{H}_\ell 
    \approx \mathbb{I} \,.
\end{equation}
That is, these discrete Hodge star operators are only approximate inverses of each other. As these duality maps are meant to be isomorphisms between the primal and dual complex, the failure of the natural Hodge star operators $\mathbb{H}_\ell$ and $\tilde{\mathbb{H}}_{n-\ell}$ to be exact inverses of each other causes one generally to avoid using both operators, but instead choose one and define the other to be its inverse. In one dimension, on a periodic domain with a uniform grid at lowest order, one finds the stencils
\begin{equation}
    \mathbbm{h}_0
    =
    \tilde{\mathbbm{h}}_0
    =
    \frac{\Delta x}{8}
    \begin{pmatrix}
        1 & 6 & 1
    \end{pmatrix}
    \,,
    \quad \text{and} \quad
    \mathbbm{h}_1
    =
    \tilde{\mathbbm{h}}_1
    =
    \Delta x^{-1}
    \mathbb{I} \,.
\end{equation}
In order to use a diagonal duality matrix, one may use $\tilde{\mathbbm{h}}_1^{-1} = \Delta x \mathbb{I}$ as the duality map on $0$-forms and $\mathbbm{h}_1 = \Delta x^{-1} \mathbb{I}$ as the duality map on $1$-forms. These diagonal matrices are the one-dimensional duality structures used in the numerical examples in this work. 

In higher dimensions on tensor-product domains, the duality matrices are simply Kronecker products of the one-dimensional matrices. In three-dimensions
\begin{multline}
    \mathbb{H}_0
    =
    \mathbbm{h}_0^x \otimes \mathbbm{h}_0^y \otimes \mathbbm{h}_0^z \,,
    \quad
    \mathbb{H}_1
    =
    \begin{pmatrix}
        \mathbbm{h}_1^x \otimes \mathbbm{h}_0^y \otimes \mathbbm{h}_0^z & 0 & 0 \\
        0 & \mathbbm{h}_0^x \otimes \mathbbm{h}_1^y \otimes \mathbbm{h}_0^z & 0 \\
        0 & 0 & \mathbbm{h}_0^x \otimes \mathbbm{h}_0^y \otimes \mathbbm{h}_1^z 
    \end{pmatrix} \,, \\
    \quad 
    \mathbb{H}_2
    =
    \begin{pmatrix}
        \mathbbm{h}_0^x \otimes \mathbbm{h}_1^y \otimes \mathbbm{h}_1^z & 0 & 0 \\
        0 & \mathbbm{h}_1^x \otimes \mathbbm{h}_0^y \otimes \mathbbm{h}_1^z & 0 \\
        0 & 0 & \mathbbm{h}_1^x \otimes \mathbbm{h}_1^y \otimes \mathbbm{h}_0^z 
    \end{pmatrix} \,,
    \quad \text{and} \quad
    \mathbb{H}_3 = \mathbbm{h}_1^x \otimes \mathbbm{h}_1^y \otimes \mathbbm{h}_1^z \,,
\end{multline}
where $\mathbbm{h}_0^\alpha$ and $\mathbbm{h}_1^\alpha$ denote the one-dimensional duality matrices in each coordinate direction. The two dimensional case is analogous.

\section{On the stability of the nonlinear time-step} \label{appendix:nonlinear_stability}

Suppose one has an update map for $\bm{\uprho}$ given by
\begin{equation}
    \bm{\uprho} \mapsto  \bm{\uprho} + \Delta t \mathbb{A} \bmsf{p}_k(\bm{\uprho}) \,,
\end{equation}
where $\mathbb{A}$ is some linear map, and $\bmsf{p}_k(\bm{\uprho})$ is a nonlinear operator. Suppose $\bm{\uprho} - \bm{\uprho}_* = \epsilon_\rho$. Linearizing about the current state of the system, one finds
\begin{equation}
    \epsilon_\rho
    \mapsto
    \epsilon_\rho
    +
    \Delta t \mathbb{A}
    (\bmsf{p}(\bm{\uprho})
    -
    \bmsf{p}(\bm{\uprho}_*)) 
    =
    (1 + \Delta t \mathbb{A} D \bmsf{p}(\bm{\uprho})) \epsilon_\rho + O(\epsilon_\rho^2) \,.
\end{equation}
In this work, one has (modulo rescaling of the density field) that
\begin{equation}
    \bmsf{p}(\bm{\uprho}) = \frac{1 - \sqrt{1 - 4 k \bm{\uprho}}}{2 k} 
    \implies
    D \bmsf{p}(\bm{\uprho})
    =
    \frac{1}{\sqrt{1- 4 k \bm{\uprho}}} \,.
\end{equation}
For the linearized time-step, Von Neumann stability analysis implies that the time-step is stable if
\begin{equation}
    \frac{\Delta t \| \mathbb{A} \|}{\sqrt{1 - 4 k \rho_*}} < 1 \,,
\end{equation}
where $\rho_* = \max_i \uprho_i$. Hence, if $k \ll 1$, and $\rho_* = O(1)$, it is sufficient to only consider the stability conditions for the time-step in the limit $k \to 0$. While these arguments are only heuristic in nature, it was found empirically in this paper's numerical experiments that it was sufficient to pick a time-step which ensures stability of the linearized dynamics. 

\end{document}